\documentclass[a4paper]{amsart}
\usepackage[utf8]{inputenc}
\usepackage[english]{babel}

\usepackage{amscd,amssymb}
\usepackage{amsmath}
\usepackage{amsthm}
\usepackage{graphicx}
\usepackage{fancybox}
\usepackage{epic,eepic}
\usepackage{amstext}
\usepackage{mathtools}
\usepackage{esint}
\usepackage[square,numbers]{natbib}
\usepackage{booktabs}
\usepackage[hide links]{hyperref}
\usepackage{comment}
\usepackage{mathtools}
\usepackage{tikz,pgfplots}
\usepackage{subcaption}

\usepackage[noabbrev, capitalise, nameinlink]{cleveref}
\crefname{equation}{}{}
\crefname{assumption}{Assumption}{Assumptions}
\crefformat{equation}{\textup{#2(#1)#3}}

\newtheorem{theorem}{Theorem}[section]
\newtheorem{corollary}[theorem]{Corollary}
\newtheorem{lemma}[theorem]{Lemma}

\theoremstyle{definition}

\theoremstyle{remark}

\numberwithin{theorem}{section}
\numberwithin{equation}{section}
\numberwithin{figure}{section}

\def\TH{\mathcal T_H}

\def\diam{\operatorname{diam}}

\def\div{\nabla\cdot}
\def\with{\,:\,}
\def\dx{\,\mathrm{d}x}
\def\ds{\,\mathrm{d}\sigma}
\def\tint{\begingroup\textstyle \int\endgroup}
\def\Nb{\mathsf{N}}

\numberwithin{equation}{section}
\numberwithin{theorem}{section}

\AtBeginDocument{%
	\def\MR#1{}
}

\title[A LOD Method for Heterogeneous Stokes Problems]{A Localized Orthogonal Decomposition Method for Heterogeneous Stokes Problems}
\author[M.~Hauck, A. Lozinski]{Moritz Hauck$^*$, Alexei Lozinski$^\dagger$}
\address{${}^+$ Interdisciplinary Center for Scientific Computing (IWR), Heidelberg University, 69047 Heidelberg, Germany}
\email{moritz.hauck@iwr.uni-heidelberg.de}
\address{${}^{\dagger}$ Université de Franche-Comté, CNRS, LmB, 25000 Besançon, France}
\email{alexei.lozinski@univ-fcomte.fr}

\begin{document}
	
\begin{abstract}
In this paper, we propose a multiscale method for heterogeneous Stokes problems. The method is based on the Localized Orthogonal Decomposition (LOD) methodology and has approximation properties independent of the regularity of the coefficients. We apply the LOD to an appropriate reformulation of the Stokes problem, which allows us to construct exponentially decaying basis functions for the velocity approximation while using a piecewise constant pressure approximation. The exponential decay motivates a localization of the basis computation, which is essential for the practical realization of the method. We perform a rigorous a priori error analysis and prove optimal convergence rates for the velocity approximation and a post-processed pressure approximation, provided that the supports of the basis functions are logarithmically increased with the desired accuracy. Numerical experiments support the theoretical results of this paper.
\end{abstract}


	
\maketitle

\section{Introduction}
This paper considers a heterogeneous Stokes problem posed in a bounded Lipschitz polytope $\Omega \subset \mathbb{R}^n$, $n \in \{2,3\}$, which we assume to be of unit size. Given an  external force~$f$, it seeks a  velocity $u$ and a pressure $p$ such that
\begin{equation}
	\label{pbStokes}
	\left\{
	\begin{aligned}
		- \div (\nu \nabla u) + \sigma
		u + \nabla p = f & \qquad \text{in } \Omega,\\
		\div u = 0 & \qquad \text{in } \Omega,\\
		u = 0 & \qquad \text{on } \partial \Omega
	\end{aligned}\right.
\end{equation}
holds, where the coefficients $\nu$ and~$\sigma$ encode the heterogeneity of the problem.~They may possibly be rough and involve oscillations on multiple non-separated length scales. 
Such problems may arise, for example, as part of more complex coupled problems where the viscosity is an unknown. A specific example is magma modeling, where the temperature-dependent viscosity can exhibit large spatial variations, cf.~\cite{Gutirrez2010}. Problems like \cref{pbStokes} can also be used as approximations to (slow) flow problems around numerous obstacles of possibly small diameter, cf.~\cite{angot99}. In this case, one sets $\nu$ as the physical viscosity and $\sigma$ to zero in the fluid domain, while assigning large values to these coefficients inside the obstacles. 

The numerical treatment of such heterogeneous problems with classical finite element methods (FEMs) suffers from suboptimal approximation rates and preasymptotic effects on meshes that do not resolve the coefficients. Since globally resolving all microscopic details of the coefficients may not be computationally feasible, we aim to construct a numerical method with reasonable errors already on coarse meshes. 
For diffusion-type problems, there is a whole zoo of multiscale methods that use coarse problem-adapted ansatz spaces.
Two classes of methods may be distinguished: First, methods that exploit structural properties of the coefficients, such as periodicity and scale separation, to construct the problem-adapted basis functions. Their computational cost differs from that of classical FEMs on the same mesh only by the cost of solving a fixed number of local problems.  This class includes the Heterogeneous Multiscale Method~\cite{EE03}, the Two-Scale FEM \cite{MaS02}, and the Multiscale FEM~\cite{HoW97}.
In contrast, methods in the second class achieve optimal approximation orders under minimal structural assumptions on the coefficients. This is achieved at the expense of a moderate computational overhead compared to classical FEMs. This overhead manifests itself in an enlarged support of the basis functions or an increased number of basis functions per mesh entity.  Prominent methods for diffusion-type multiscale problems include the Generalized Multiscale FEM~\cite{EfeGH13,ChuEL18b}, the Multiscale Spectral Generalized FEM  (MS-GFEM) \cite{BabL11,Ma22}, Adaptive Local Bases~\cite{GraGS12}, the \text{(Super-)} Localized Orthogonal Decomposition (LOD) method~\cite{MalP14,HenP13,HaPe21b}, or Gamblets~\cite{Owh17}; see also the review article~\cite{AltHP21}. 
We mention that the LOD  and the MS-GFEM have been applied to Darcy-type problems, cf.~\cite{Mlqvist2016,Alber2024}, which is also of mixed form but different in nature than the Stokes problem; see also \cite{LARSON2009}.

Some of the methods of the first class have been successfully generalized to  heterogeneous Stokes problems. For example, a homogenization-based method has been proposed in \cite{Brown2013,Brown12013} for slowly varying perforated media. As for the Multiscale FEM, its Crouzeix--Raviart version, first proposed in \cite{LeBris2014}, has been applied to the Stokes problem in perforated domains in \cite{Muljadi2015,Jankowiak2024,Feng22,balaziPhD}. The weak notion of continuity of the Crouzeix--Raviart method and an appropriate choice of the problem-adapted approximation space make it flexible enough to cope with the divergence-free constraint of the Stokes problem. On the contrary, we are not aware of any multiscale method for heterogeneous Stokes problems that works under the minimal structural assumptions on the coefficients. 

This open question is addressed in this paper by adapting the LOD methodology to heterogeneous Stokes problems. The basic idea of the LOD is to decompose the solution space into a fine-scale space and its orthogonal complement with respect to the energy inner product induced by the considered problem. By choosing the fine-scale space to consist of functions that average out on coarse scales, one obtains a finite-dimensional mesh-based complement space that is adapted to the problem at hand and has uniform approximation properties under minimal structural assumptions on the coefficients. It possesses exponentially decaying basis functions whose computation can thus be localized to subdomains, resulting in a practically feasible method.  For Stokes problems, the divergence-free constraint poses a major challenge to LOD-type methods, since directly incorporating the constraint can lead to ill-posed problems or slowly decaying basis functions. To overcome this problem, we reformulate the Stokes problem using the space of $(H^1_0(\Omega))^n$-functions, whose divergence is piecewise constant with respect to some coarse mesh, as the solution space for the velocity. The divergence-free velocity is then recovered using a piecewise constant Lagrange multiplier defined on the same mesh. Inspired by the Crouzeix--Raviart Multiscale FEM mentioned above, we then choose the fine-scale space for the velocity as the functions whose averages vanish on all faces of the coarse mesh. We prove that the finite-dimensional orthogonal complement possesses exponentially decaying basis functions, which paves the way to the construction of an LOD-type multiscale method. The use of face averages to define the fine-scale space is novel in the LOD context and allows us to construct the velocity approximation space independent of the pressure. The  approximation of the resulting LOD method is exactly divergence-free and thus pressure robust in the sense of \cite{John2017}. 
We also perform an \textit{a priori} error analysis of the proposed method and prove optimal orders of convergence as the mesh size is decreased, provided that the support of the basis functions is allowed to increase logarithmically with the desired accuracy. More specifically, for $L^2$-right-hand sides we prove first- and second-order convergence for the $L^2$- and $H^1$-errors of the velocity approximation, respectively, and first-order convergence for a post-processed pressure approximation. If the right-hand side is $H^1$-regular, we can squeeze out an additional order of convergence for the velocity approximation.  We emphasize that only minimal structural assumptions on the coefficients are necessary for this error analysis. 

The paper is organized as follows: In \cref{sec:modelproblem} we introduce the model problem and a prototypical multiscale method is presented in \cref{sec:protmethod}. We prove the exponential decay of the prototypical basis functions and localize the basis computation in \cref{sec:expdec}. A practical multiscale method is then presented in \cref{sec:pracmethod}. Numerical experiments supporting our theoretical results are given in \cref{sec:numexp}.

	\section{Model problem}
	\label{sec:modelproblem}
This section introduces the weak formulation of the heterogeneous Stokes problem, along with classical results guaranteeing its well-posedness. The weak formulation is based on the Sobolev space 
$V \coloneqq (H^1_0(\Omega))^n$ endowed with homogeneous Dirichlet boundary conditions on $\partial \Omega$ and the space $M \coloneqq \{q \in L^2(\Omega) \with \int_\Omega q\dx = 0\}$ satisfying an integral-mean-zero constraint. In the following, we will always assume that there exist constants $\nu_\mathrm{min}, \nu_\mathrm{max}$ and $\sigma_\mathrm{min}, \sigma_\mathrm{max}$ such that 
\begin{equation}
	\label{eq:unifboundcoeff}
	0<\nu_{\min}\le\nu\le\nu_{\max}<\infty,\qquad 0\leq \sigma_{\min}\le\sigma\le\sigma_{\max}<\infty
\end{equation}
holds almost everywhere in $\Omega$. Denoting by  $(\cdot,\cdot)_\Omega$ the $L^2(\Omega)$-inner product, the problem's bilinear forms $a\colon V\times V \to \mathbb R$ and $b\colon V\times M \to \mathbb R$ are   defined as
\begin{equation*}
	a(u,v) \coloneqq (\nu \nabla u ,\nabla v)_{\Omega}  + (\sigma u, v)_{\Omega},\qquad b(u,q) \coloneqq -(q,\div u)_\Omega.
\end{equation*}

Given a source term $f \in L^2(\Omega)$, the weak formulation of the considered heterogeneous Stokes problem seeks a pair $(u,p) \in V \times M$ such that 
\begin{subequations}
	\label{eq:weakstokes}
	\begin{align}
		\qquad\qquad \qquad \qquad &a(u,v)& +\quad &b(v,p) &=\quad &(f,v)_\Omega,&\qquad \qquad \qquad\qquad\label{eq:weakstokes1}\\
		\qquad \qquad \qquad&b(u,q)&   &       &=\quad &0&\qquad \qquad \qquad\label{eq:weakstokes2}
	\end{align}
\end{subequations}
holds for all $(v,q) \in V \times M$.

Using the uniform coefficient bounds~\cref{eq:unifboundcoeff} one can show that the bilinear form~$a$ is continuous and coercive, i.e., there exist constants $c_a, C_a>0$ such that
\begin{equation}
	\label{eq:a}
	 |a(v,v)|\geq c_a \|\nabla v\|_\Omega^2,\qquad|a(u,v)| \leq C_a \|\nabla u\|_\Omega\|\nabla v\|_\Omega
\end{equation}
holds for all $u,v \in V$, where $\|\cdot\|_\Omega^2 \coloneqq (\cdot,\cdot)_\Omega$ denotes the $L^2(\Omega)$-norm. Note that by the Poincaré--Friedrichs inequality, the seminorm $\|\nabla \cdot \|_\Omega$ is equivalent to the full $(H^1(\Omega))^n$-norm. The constants in \cref{eq:a} can be specified as $c_a = \nu_\mathrm{min}$ and $C_a = \nu_\mathrm{max} + C_\mathrm{PF}^2\sigma_\mathrm{max}$, where $C_\mathrm{PF}>0$ denotes the Poincaré--Friedrichs constant.  

To establish the well-posedness of problem \cref{eq:weakstokes}, we  need a compatibility condition between the spaces $V$ and $M$, expressed as the inf--sup condition
\begin{equation}
	\label{eq:infsup}
	\adjustlimits \inf_{q \in M}\sup_{v \in V} \frac{|b(v,q)|}{\|\nabla v\|_\Omega\|q\|_\Omega}\geq  c_b,
\end{equation}
where $c_b>0$ is typically called the inf--sup constant. This condition is classical and it is typically proved using the so-called Ladyzhenskaya lemma, cf.~\cite{ladyzhenskaia1963}. It states that for any $q \in M$ there exists $v \in V$ such that
\begin{equation}
	\label{eq:ladlemmaOmega}
	\div v = q,\qquad \|\nabla v\|_\Omega \leq C_\mathrm{L} \|q\|_\Omega,
\end{equation}
which directly implies the inf--sup stability with inf--sup constant $c_b = C_\mathrm{L}^{-1}$. After establishing conditions \cref{eq:a,eq:infsup}, the well-posedness of weak formulation~\cref{eq:weakstokes} can be concluded using classical inf--sup theory; see, e.g.,~\cite{BoffiBrezziFortin2013}.

	\section{Prototypical multiscale method}
	\label{sec:protmethod}
	
	This section introduces a prototypical multiscale method that achieves optimal order approximations without preasymptotic effects under minimal structural assumptions on the coefficients. 
To this end, we introduce a geometrically conforming, quasi-uniform, and shape-regular hierarchy of simplicial\footnote{One can also use general polygonal/polyhedral meshes. We assume that the meshes are simplicial only to simplify the presentation.  The method itself can be applied to more general meshes in a straightforward manner. An extension to meshes with curved elements  is also possible.}  meshes $\{\TH\}_{H}$. 
Each mesh is a finite subdivision of the closure of $\Omega$, into closed elements~$K$, which are $n$-dimensional simplices. The parameter~${H>0}$ denotes the  mesh size and is defined as the maximum diameter of the elements in $\TH$, i.e., $H \coloneqq \max_{K \in \TH} \diam(K)$. We further denote the space of $\TH$-piecewise constant functions by~$\mathbb P^0(\TH)$ and  write $\Pi_H\colon L^2(\Omega)\to \mathbb P^0(\TH)$ for the corresponding $L^2$-orthogonal projection. The set of all faces of the mesh~$\TH$ is denoted by $\mathcal F_H$ and the subset of interior faces by $\mathcal F_H^i$.

For the construction of the prototypical multiscale method we will use an equivalent reformulation of problem~\cref{eq:weakstokes}. This reformulation is based on the spaces
	\begin{equation}
		\label{eq:defZ}
		Z \coloneqq \big\{ v \in V \with \div v \in \mathbb P^0(\TH)\big\},\qquad M_H \coloneqq M \cap \mathbb P^0(\TH),
\end{equation}
where the space $Z$ partially integrates the divergence-free constraint into the velocity space. Thus, the smaller space $M_H$ is sufficient to enforce that the velocity is divergence-free. The reformulation seeks  $(u,p_H) \in
	Z \times M_H$ such that
	\begin{subequations}
		\label{eq:reformulation}
			\begin{align}
				\quad\qquad \qquad \qquad &a(u,v)& +\quad &b(v,p_H) &=\quad &(f,v)_\Omega,&\qquad \qquad \qquad\qquad\label{eq:reformulation1}\\
			\quad \qquad \qquad&b(u,q_H)&   &       &=\quad &0.&\qquad \qquad \qquad
		\end{align}
	\end{subequations}
holds for all $(v,q_H) \in Z\times M_H$. 
To prove the well-posedness of this reformulated problem we verify the corresponding inf--sup condition. This inf--sup condition can be shown to hold with the constant $c_b$ from~\cref{eq:infsup}, using again the Ladyzhenskaya lemma, cf.~\cref{eq:ladlemmaOmega}. It is easy to verify that the first component of the solution to the reformulated problem coincides with the velocity $u$ from \cref{eq:weakstokes} and that for the second component we have $p_H = \Pi_H p$, where $p$ is the pressure from \cref{eq:weakstokes}.

Following the construction of the LOD for diffusion-type problems, cf.~\cite{MalP14, MalP20}, we consider a decomposition of the space $Z$ into the direct sum of two subspaces. The first one, typically referred to as fine-scale space, contains functions that average out on coarse scales and is defined as
	\begin{equation}
		W \coloneqq \big\{ v \in Z\with\tint_F v\ds = 0,\; F \in \mathcal F_H^i \big\}.
	\end{equation}
The second subspace is finite-dimensional and will serve as the approximation space of the prototypical LOD method. It is defined as the orthogonal complement of $W$ with respect to the energy inner product $a$, i.e.,
	\begin{equation}
		\label{eq:Zms}
		\tilde Z_H \coloneqq \big\{ u \in Z \with a (u, v) = 0,\; v \in W \big\}.
	\end{equation}
Note that, since $\tilde Z_H$ is constructed as the orthogonal complement of $W$ with respect to the problem-dependent inner product $a$, it contains problem-specific information that allows reliable approximations even at coarse scales. The use of tildes in the notation of functions and spaces is intended to emphasize that they are adapted to the problem at hand.  The following lemma constructs a basis of the space $\tilde Z_H$. 

\begin{lemma}[Prototypical basis]\label{le:protbasis}The space $\tilde Z_H$ has dimension $N \coloneqq n \cdot \# \mathcal F_H^i$ with~$\# \mathcal F_H^i$ denoting the number of interior faces, and a basis of it is given by
	\begin{equation}
		\label{eq:LODprotbasis}
		\big\{\tilde \varphi_{F,j}\with F \in \mathcal F_H^i,\; j =1,\dots,n\big\}
	\end{equation}
	with $\tilde \varphi_{E,j}$ defined for all $F \in \mathcal F_H^i$ and $j =1,\dots,n$ as the unique solutions to: seek $(\tilde \varphi_{F,j},\xi_{F,j},\lambda) \in V\times Q\times \mathbb R^N$ with $Q \coloneqq \{ q \in M \with \int_K q
	\dx= 0 ,\; K \in \mathcal{T}_H \}$ such that
	\begin{subequations}
		\label{pbphiE} 
		\begin{align}
			&\quad\qquad \qquad a (\tilde \varphi_{F,j}, v)& +\quad  &b(v,\xi_{F,j})          & +\quad  &c(v,\lambda) & =\quad  &0,\qquad \quad\quad &\label{eq:LODbasis1}\\
			& \quad\qquad \qquad  b(\tilde \varphi_{F,j},\chi)                &   &&   &             & =\quad  &0,\qquad\quad\qquad &\label{eq:LODbasis2}\\
			&\quad\qquad\qquad c(\tilde \varphi_{F,j},\mu) &   &                  &   &             & =\quad  &\mu_{F, j}\qquad \quad\quad\label{eq:LODbasis3}&
		\end{align}
	\end{subequations}
	holds for all $(v,\chi,\mu) \in V\times Q\times \mathbb R^N$. Here, we label the entries of a vector $\mu\in\mathbb R^N$ using face-index pairs as $\mu=\{\mu_{E,k}\with E \in \mathcal F_H^i,\; k =1,\dots,n\big\}$, and the bilinear form $c \colon V \times \mathbb R^N\to \mathbb R$ is defined as
	\begin{equation*}
		c(v,\mu) \coloneqq \sum_{E \in \mathcal F_H^i} \sum_{k = 1}^n \mu_{E, k}  \int_E v
		\cdot e_k \ds,
	\end{equation*}
	where $\{e_k\with k=1,\dots,n\}$ denotes the canonical basis of $\mathbb R^n$.
\end{lemma}

Before proving this lemma, we introduce some technical tools that will be used not only in its proof, but also on several other occasions in the remainder of this paper.  First, we note that the Ladyzhenskaya lemma stated in~\cref{eq:ladlemmaOmega} for the whole domain $\Omega$ also holds locally on elements $K \in \TH$. More precisely, there exists a constant~${C_\mathrm{L}'>0}$ such that, for every element $K \in \TH$ and all functions $q \in L^2(K)$ with $\int_Kq\dx =0$, there exists $v \in (H^1_0(K))^n$ such that 
\begin{equation}
	\label{eq:ladlemma}
	\div v = q,\qquad \|\nabla v\|_K \leq C_\mathrm{L}' \|q\|_K
\end{equation}
holds, where the constant depends only on the shape regularity of the mesh.
The latter statement can be inferred, for example,  from~\cite{bernardi16}, where the shape-dependence of the Ladyzhenskaya constant is investigated. The locally supported functions $v$ from \cref{eq:ladlemma}, hereafter referred to as element bubble functions, will be frequently used to estimate the bilinear form~$b$. 

Furthermore, for estimating the bilinear form $c$, we introduce face bubble functions denoted by $\{b_{F,j}\with F \in \mathcal F_H^i,\; j = 1,\dots,n\}$. Each bubble function is locally supported with~$b_{F,j} \in (H^1_0(\omega_F))^n$, where $\omega_F$ is the union of the two elements that share face~$F$. The edge bubbles are chosen such that ${\tint_E b_{F,j}\cdot e_k \ds = \delta_{EF}\delta_{jk}}$ holds for all $E,F \in \mathcal F_H^i$ and $j,k = 1,\dots,n$, and the two stability estimates
\begin{equation}
	\label{eq:bubbleest}
	\|b_{F,j}\|_{\omega_F} \leq C_\mathrm{b} H^{-n/2+1},\qquad \|\nabla b_{F,j}\|_{\omega_F} \leq C_\mathrm{b} H^{-n/2}
\end{equation}
 are satisfied for a constant $C_\mathrm{b}>0$ independent of $H$. Additionally, we demand that the divergence of the edge bubbles is piecewise constant, i.e., ${b_{F,j} \in Z}$. 
Note that bubbles with these properties can be constructed using classical edge bubbles, which we denote by $\psi_F$ below, cf.~\cite{Verfrth2019}. They satisfy $\tint_E \psi_{F}\ds = \delta_{EF}$ for any ${E,F \in \mathcal F_H^i}$. We then define $b_{F,j}\coloneqq \psi_Fe_j+v_{F,j}$, where $v_{F,j}$ is supported on $\omega_F$ and defined locally on $K \subset \omega_F$ as the $(H^1_0(K))^n$-function satisfying $\div v_{F, j}=-\div(\psi_Fe_j)+n_{F,K}\cdot e_j$, where $n_{F,K}$ denotes the unit normal on $F$ outward of~$K$. The function~$v_{F,j}$ exist thanks to~(\ref{eq:ladlemma}) and satisfies $\|\nabla v_{F,j}\|_K \leq C_\mathrm{L}' \|\nabla \psi_F\|_K$.

Now let us proceed with the proof of the lemma above.
\begin{proof}[Proof of \cref{le:protbasis}]
	We first prove the well-posedness of problem~\cref{pbphiE}. To this end, we verify the inf--sup stability of the  bilinear form $d\colon V \times (Q\times \mathbb R^N)$ defined as $d(v,(\chi,\mu)) \coloneqq b(v,\chi) + c(v,\mu)$, i.e., there exists a constant $c_d>0$ such that
\begin{equation}\label{eq:infsupdoublespp}
	\adjustlimits \inf_{(\chi,\mu) \in Q \times \mathbb R^N}\sup_{v \in V} \frac{|d(v,(\chi,\mu))|}{\|\nabla v\|_\Omega \|(\chi,\mu)\|}\geq c_d
\end{equation}
holds, where $\|(\chi,\mu)\|^2 \coloneqq \|\chi\|_\Omega^2+|\mu|^2$ with $|\cdot|$ denoting the Euclidean norm.

Given any $(\chi,\mu) \in Q\times \mathbb R^N$, we choose a function $v \in V$ defined as
\begin{align*}
	v \coloneqq \sum_{K \in \mathcal T_H} v_K + \sum_{F \in \mathcal F_H^i} \sum_{j = 1}^n \mu_{F, j} b_{F,j},
\end{align*}
where, for any $K\in\TH$, $v_K \in (H^1_0(K))^n$ satisfies $\div v_K = \chi$ locally on $K$ as well as $\|\nabla v_K\|_K \leq C_\mathrm{L}' \|\chi\|_K$, cf.~\cref{eq:ladlemma}.  By construction, it holds $|d(v,(\chi,\mu))| = \|(\chi,\mu)\|^2$ and we have for any $K\in\TH$ that
\begin{align*} \| \nabla v\|_K^2 &\leq 2 (C_\mathrm{L}^{\prime})^2 \| \chi \|_K^2 + 2 \bigg(
\sum_{F \in \mathcal{F}_H^i \cap \partial K} \sum_{j = 1}^n \mu_{F, j}^2
\bigg) \bigg( \sum_{F \in \mathcal{F}_H^i \cap \partial K} \sum_{j = 1}^n 
\| \nabla b_{F, j} \|_K^2 \bigg) \\
 &\leq 2 (C_\mathrm{L}^{\prime})^2 \| \chi \|_K^2 + 2 n (n + 1) C_{\mathrm{b}}^2
H^{- n} \bigg( \sum_{F \in \mathcal{F}_H^i \cap \partial K} \sum_{j = 1}^n
\mu_{F, j}^2 \bigg),
\end{align*}
where we used estimates~\cref{eq:bubbleest,eq:ladlemma}. Summing over all elements gives
\begin{equation*}
 \| \nabla v\|_{\Omega}^2 \leq \max \{2 (C_\mathrm{L}^{\prime})^2, 4 n (n + 1)
C_{\mathrm{b}}^2 H^{- n}\} \|(\chi, \mu)\|^2,
\end{equation*}
which proves inf--sup condition~\cref{eq:infsupdoublespp} with the constant $c_d$, defined as the reciprocal of the constant in the estimate above.

Next, we prove that~\cref{eq:LODprotbasis} constitutes a basis of the space $\tilde Z_H$. To this end, we first note that conditions \cref{eq:LODbasis1,eq:LODbasis2} imply that $\tilde \varphi_{F,j} \in \tilde Z_H$  for all $F \in \mathcal{F}_H^i$ and $j \in \{1, \ldots,n\}$. To prove the basis property, we consider an arbitrary $u \in \tilde{Z}_H$ and set $w = u - \sum_{F \in \mathcal{F}_H^i} \sum_{j = 1}^n  \int_F u \cdot e_j \ds \tilde{\varphi}_{F,	j}$. Since $w \in W$ and, at the same time, $w$ is $a$-orthogonal to $W$, it  must hold that $w = 0$. Thus, any $u \in \tilde{Z}_H$ can be written as a unique linear combination of the linearly independent functions~$\tilde \varphi_{F,j}$, which shows that they form a basis of the space  $\tilde Z_H$. This concludes the proof.
\end{proof}

Having introduced the prototypical basis functions, we can define a projection operator $\mathcal R \colon Z \to \tilde Z_H$ which preserves face integrals for any $v \in Z$ by
\begin{equation}
	\label{eq:defR}
	\mathcal R v \coloneqq \sum_{F \in \mathcal F_H^i} \sum_{j = 1}^n  
	\int_F v \cdot e_j \ds\, \tilde \varphi_{F,j}.
\end{equation}
This operator coincides with the $a$-orthogonal projection onto $\tilde Z_H$, since for any $v \in Z$ and $w \in \tilde Z_H$ it holds that $a(v-\mathcal R v,w)$, noting that $v-\mathcal R v\in W$. This orthogonality property also implies the continuity of  $\mathcal R$, i.e., for any $v \in Z$ it holds
\begin{equation}
	\label{eq:Rcont}
	\|\nabla \mathcal R v\|_\Omega \leq \sqrt{C_a/c_a} \|\nabla v\|_\Omega,
\end{equation}
where we used~\cref{eq:a}. 
Furthermore, applying the divergence theorem, using the fact that~$\mathcal R$ preserves face integrals, and noting that by the definition of $Z$ functions in this space have a piecewise constant divergence, cf.~\cref{eq:defZ}, one can show that 
\begin{equation}
	\label{eq:identityIH}
	 (\div \mathcal R v) |_K = (\div v)|_K
\end{equation}
holds for any $v \in Z$ and $K \in \TH$. 

The  desired prototypical method then seeks  $(\tilde u_H, \tilde p_H) \in \tilde Z_H \times M_H$ such that
	\begin{subequations}
	\label{LODid2}
	\begin{align}
 	\qquad \qquad \qquad &a(\tilde u_H,\tilde v_H)& +\quad &b(\tilde v_H,\tilde p_H) &=\quad &(f,\tilde v_H)_\Omega,&\qquad \qquad \qquad\label{eq:LODdivfree1}\\
 \qquad \qquad \qquad&b(\tilde u_H,\tilde q_H)&   &       &=\quad &0.&\qquad \qquad \qquad\label{eq:LODdivfree2}
		\end{align}
	\end{subequations}
holds for all $(\tilde v_H,\tilde q_H) \in \tilde Z_H\times M_H$. 
For the analysis of this method in the theorem below we need an approximation result for $\Pi_H$. Using the Poincaré inequality on convex domains, cf.~\cite{PaW60}, it hollows that for all $K \in \TH$ and $v \in H^1(K)$ it holds
\begin{equation}
	\label{eq:Poincare}
	\|v-\Pi_H v\|_K \leq \pi^{-1}H\|\nabla v\|_K.
\end{equation}

\begin{theorem}[Prototypical method]\label{thm:convergenceprot}
	The prototypical multiscale method~\cref{LODid2} is well-posed, and its solution is given by $(\tilde u_H,\tilde p_H)=(\mathcal R u,\Pi_Hp)$, where $(u,p)$ solves problem~\cref{eq:weakstokes}.
Furthermore, there exist constants $C_\mathrm{e}, C_\mathrm{e}'>0$ independent of $H$ such that for any right-hand side $f \in H^s(\Omega)$ with $s \in \{0,1\}$, we have the error~estimates
		\begin{align}
			\|u - \tilde u_H\|_\Omega &\leq C_\mathrm{e}H^{2+s} |f|_{s,\Omega},\label{eq:errestL2}\\
			\|\nabla (u - \tilde u_H)\|_\Omega &\leq C_\mathrm{e}' H^{1+s} |f|_{s,\Omega},\label{eq:errestH1}
		\end{align}
		where $|\cdot|_{s,\Omega}$ denotes the $H^s(\Omega)$-seminorm.
\end{theorem}
\begin{proof}

First, we prove the  inf--sup condition
\begin{equation}\label{eq:infsupLOD}
	\adjustlimits \inf_{\tilde q_H \in M_H}\sup_{\tilde v_H \in \tilde Z_H} \frac{|b(\tilde v_H,\tilde q_H)|}{\|\nabla \tilde v_H\|_\Omega\|\tilde q_H\|_\Omega} \geq \tilde c_b
\end{equation}
for some constant $\tilde c_b>0$, which implies the well-posedness of~\cref{LODid2}.
	Given any $\tilde q_H \in M_H$, we denote by $v \in V$ the function with $\div v = \tilde q_H$ and $\|\nabla v\|_\Omega \leq C_\mathrm{L} \|\tilde q_H\|_\Omega$, cf.~\cref{eq:ladlemmaOmega}, and define $\tilde v_H \coloneqq \mathcal R v$. One can show that $|b(\tilde v_H,\tilde q_H)| = |b(v,\tilde q_H)| = \|\tilde q_H\|_\Omega^2$ using~\cref{eq:identityIH}, and that $\|\nabla \tilde v_H\|_\Omega \leq \sqrt{C_a/c_a} \|\nabla  v\|\leq C_\mathrm{L}\sqrt{C_a/c_a}\|\tilde q_H\|_\Omega$ by the continuity of $\mathcal R$, cf.~\cref{eq:Rcont}, and the particular choice of $v$. Combining these results, inf--sup condition~\cref{eq:infsupLOD} follows directly with the constant $\tilde c_b = (C_\mathrm{L}\sqrt{C_a/c_a})^{-1}$.
	Using that $\mathcal R\colon Z \to \tilde Z_H$ is $a$-orthogonal, identity~\cref{eq:identityIH}, as well as that $(u,p_H)$ solves~\cref{eq:reformulation}, it is a straightforward observation that $(\mathcal Ru,p_H)$ solves \cref{LODid2}. The uniqueness of the solution to \cref{LODid2} then implies that $\tilde u_H = \mathcal Ru$ and $\tilde p_H = p_H$. 
	
	Second, we derive error estimates for the prototypical multiscale method. To this end, we note that the  error $e \coloneqq u-\tilde u_H$ is an element of the space~$W$. Furthermore, we have for all $K \in \TH$ that $(\div \tilde u_H)|_K = (\div u)|_K = 0$, where we apply~\cref{eq:identityIH} to~$u$ and note that $\tilde u_H = \mathcal Ru$. This implies that $\div e = \div \tilde u_H = 0$. Using the first estimate in~\cref{eq:a}, that $a(\tilde u_H,e) = 0$, \cref{eq:reformulation1} with $e$ as test function, that $b(e,p_H) = 0$, and the local Poincaré-type inequality from~\cref{lem:Poincare}, we obtain that
	\begin{equation}
		\label{eq:errestprot}
		 c_a\|\nabla e\|_\Omega^2 
		\leq a (e, e) = a (u, e) = (f, e)_{\Omega}
		\leq C_\mathrm{P}  H \|f\|_\Omega \|\nabla e\|_\Omega,
	\end{equation}
which implies the desired $H^1$-error estimate for $s = 0$.

In the case $s = 1$, we add and subtract $\Pi_H f$ in~\cref{eq:errestprot}, the $L^2$ projection of $f$ onto piecewise constants, which results in
\begin{equation}
	\label{2terms}
	c_a \|\nabla e\|_\Omega^2 \leq (f-\Pi_H f,e)_\Omega + (\Pi_H f,e)_\Omega,
\end{equation} 
where the first term on the right-hand side can be bounded using approximation result \cref{eq:Poincare} and the Poincaré-type inequality from \cref{lem:Poincare} as
 \begin{equation}\label{eq:est1}
 	 (f - \Pi_H f, e)_{\Omega} \leq \|f - \Pi_H f\|_{\Omega} \|e\|_{\Omega}
 	\leq \pi^{-1}C_{\mathrm{P}} H^2 \| \nabla f\|_{\Omega}  \| \nabla e\|_{\Omega}.
 \end{equation}
 Abbreviating $f_K \coloneqq (\Pi_H f)|_K$, the second term on the right-hand side of \cref{2terms} can be locally rewritten using the divergence theorem and that $\div e = 0$ as
\begin{equation}\label{trickyCR}
	 (\Pi_H f, e)_K = \int_K \nabla (f_K \cdot x) \cdot e\dx  =
	\int_{\partial K} (f_K \cdot x)( e \cdot n) \ds.
\end{equation}
Summing over all elements using that $\int_F e\cdot n \ds  = 0$ for all $F \in \mathcal F_H^i$ then gives 
\begin{align}
	\label{eq:refPifH}
	 (\Pi_H f, e)_{\Omega} &=
	   \sum_{F \in
		\mathcal{F}_H^{i}} \int_F \big([\Pi_H f]_F \cdot (x - x_F)\big) (e \cdot n_F) \ds,
\end{align}
where  we denote by $x_F$ the
barycenter of $F$ and by $[\cdot]_F$ the jump of a function across face~$F$ chosen consistently with the fixed normal $n_F$ associated to~$F$. 

To estimate the term on the right-hand side of the previous equation, we combine the approximation result of  \cref{eq:Poincare} with   \cref{lem:trace}, which gives for any $F \in \mathcal F_H^i$ and $K \in \TH$ with $\partial K \supset F$ that $\|v-v_K\|_F^2 \leq CH\|\nabla v\|_K^2$ holds for any $v\in H^1(K)$ with $v_K \coloneqq (\Pi_Hv)|_K$ and the constant $C \coloneqq C_\mathrm{T}(1+\pi^{-1})\pi^{-1}$. Applying this estimate to $e$ and $f$, noting that $ \|e\|_F^2 \leq \|e-e_K\|_F^2 $ and  $[\Pi_H f]_F=-[f-\Pi_H f]_F$ yields that
\begin{equation*}\label{2traceIneqs}
	\|e\|_F^2
	\leq \tfrac12CH \| \nabla e\|_{\omega_F}^2,\qquad \| [\Pi_H f]_F^2
	\|_F^2 \leq 2C H \| \nabla f\|_{\omega_F}^2.
\end{equation*}
Combining \cref{eq:refPifH} with the latter bounds and using the Cauchy--Schwarz inequality, we obtain for the second term on the right-hand side of \cref{2terms} that
\begin{equation}\label{eq:est2}
	 (\Pi_H f, e)_{\Omega} 	\leq 2(n+1)C H^2 \|	\nabla f\|_{\Omega} \| \nabla e\|_{\Omega}.
\end{equation}
Inserting estimates \cref{eq:est1,eq:est2} into \cref{2terms} yields the desired $H^1$-error estimate in the case that $s = 1$. Therefore, taking the maximum of the constants yields  estimate~\cref{eq:errestH1} with $C_\mathrm{e}' \coloneqq \max\{c_a^{-1}C_\mathrm{P}, c_a^{-1}(\pi^{-1}C_\mathrm{P} + 2(n+1)C)\}$.

The $L^2$-error estimate \cref{eq:errestL2} with constant $C_\mathrm{e} \coloneqq C_\mathrm{P}C_\mathrm{e}'$ can be proved by applying \cref{lem:Poincare} once again. This completes the proof.
\end{proof}

We emphasize that squeezing out an additional order of convergence for $H^1$-regular right-hand sides only works in the context of Stokes problems and not for diffusion-type problems. This observation can also be verified numerically. 

\section{Exponential decay and localization}
	\label{sec:expdec}
	
We emphasize that the prototypical LOD basis functions defined in \cref{eq:LODprotbasis} are globally supported. Therefore, their computation would require the solution to  global problems, which we consider infeasible. In this section, we show that the prototypical LOD basis functions decay exponentially, which motivates their approximation by locally computable counterparts. A practical multiscale method based on such local approximations is presented in \cref{sec:pracmethod}.
	
	To quantify the decay of the basis functions, we introduce the notion of patches with respect to the coarse mesh $\TH$. Given an oversampling parameter $\ell \in \mathbb N$, where~$\mathbb N$ denotes the (positive) natural numbers, we define the patch of order~$\ell$ around a union of elements $S \subset \Omega$ recursively by
	\begin{equation}
		\label{eq:patch}
		\Nb^\ell(S) \coloneqq  \Nb^1(\Nb^{\ell-1}(S)),\quad \ell \geq 2,\qquad\Nb^1(S) \coloneqq  \bigcup \bigl\{K \in \TH\,\colon\, \overline{S} \,\cap\, \overline{K}\neq \emptyset\bigl\},
	\end{equation}
	and set $\Nb(S) \coloneqq \Nb^1(S)$. The notion of patches can also be extended to faces by defining, for any $F \in \mathcal F_H$, the patch of order $\ell$ around $F$ by $\Nb^\ell(F) = \Nb^{\ell}(\omega_F)$, where $\omega_F$ denotes the union of the two elements sharing the face $F$.
The following theorem proves the exponential decay of the prototypical LOD basis functions.

\begin{theorem}[Exponential decay]\label{thm:dec}
	There exists a constant $c>0$ independent of $H$, $\ell$, $F$, and $j$ such that for all $F \in \mathcal F_H^i$, $j\in \{1,\dots,n \}$, and $\ell \in \mathbb N$, it holds that 
	\begin{equation}
		\label{eq:dec}
		\|\nabla \tilde \varphi_{F,j}\|_{\Omega\setminus \Nb^\ell(F)} \leq \exp(-c \ell)\|\nabla \tilde \varphi_{F,j}\|_\Omega.
	\end{equation}
\end{theorem}

\begin{proof}
	In the following, we will use the abbreviation $\tilde \varphi \coloneqq \tilde \varphi_{F,j}$ and consider a cut-off function $\eta \in W^{1,\infty}(\Omega)$ with the properties
\begin{align}
	\label{eq:eta}
		\left\{ \begin{aligned}
			\eta & \equiv 0 & \quad & \text{in } \Nb^{\ell - 1} (F),\\
			\eta & \equiv 1 & \quad & \text{in } \Omega \setminus \Nb^{\ell}
			(F),\\
			0 & \leq \eta \leq 1 & \quad & \text{in } R \coloneqq \Nb^{\ell} (F)
			\setminus \Nb^{\ell - 1} (F)
		\end{aligned} \right.  
\end{align}
and the bound
\begin{equation}
	\label{eq:boundeta}
\sup_{x \in \Omega}|\nabla \eta|\leq C_\eta H^{-1},
\end{equation}
where $C_\eta>0$ is a constant independent of $H$.

Using $\eta \tilde{\varphi}$ as a test function in \cref{eq:LODbasis1} gives that
\begin{equation*}
 a (\tilde{\varphi} , \eta \tilde{\varphi}) =  - b(\eta \tilde{\varphi}, \xi) - c (\eta \tilde{\varphi}, \lambda),
\end{equation*}
which, using subscripts to denote the restrictions of the bilinear forms $a,b,$ and $c$ to subdomains, can be rewritten as
\begin{equation}
	\label{eq:identsubd}
	 a_{\Omega \setminus \Nb^{\ell} (F)} (\tilde{\varphi}, \tilde{\varphi}) = - a_{R} (\tilde{\varphi} , \eta \tilde{\varphi})  -
	b_{R}  (\eta \tilde{\varphi}, \xi) - c_{R} (\eta \tilde{\varphi}, \lambda).
\end{equation}
Here, we dropped the subscript of $\xi$ and used  that $\operatorname{supp}(\eta \tilde \varphi) \subset (\Omega \setminus \Nb^{\ell} (F)) \cup R$, that the divergence of $\eta \tilde \varphi$ is piecewise constant on $\Omega \setminus R$, and that $\int_E \eta \tilde \varphi\ds =0$ for all faces $E$  not inside $R$. Using coefficient bound~\cref{eq:unifboundcoeff}, we can estimate \cref{eq:identsubd} as
\begin{align}
	\label{eq:decaymanipulations}
	\begin{split}
			c_a \|\nabla \tilde \varphi\|_{\Omega \setminus \Nb^{\ell}(F)}^2 & \leq - a_{R} (\tilde{\varphi} , \eta \tilde{\varphi})  -
     b_{R}  (\eta \tilde{\varphi}, \xi) - c_{R} (\eta \tilde{\varphi}, \lambda)\\
		&\eqqcolon \Xi_1 + \Xi_2 + \Xi_3.
	\end{split}
\end{align}

To estimate $\Xi_1$, we again use coefficient bound~\cref{eq:unifboundcoeff}, the bound on $\eta$ from~\cref{eq:boundeta}, and the Poincaré-type inequality from \cref{lem:Poincare} to get that
\begin{equation}
	\label{eq:corrected}
	\Xi_1 \leq  C_a'\| \nabla \tilde{\varphi} \|_{R} \| \nabla (\eta
	\tilde{\varphi}) \|_{R} \leq C_a'(1 + C_{\eta} C_\mathrm{P}) \| \nabla
	\tilde{\varphi} \|_{R}^2
\end{equation}
with the constant $C_a' \coloneqq (\nu_\mathrm{max} + C_\mathrm{P}^2\sigma_\mathrm{max})$. Above we used that, for any element  $K \subset R$, there exists a face $E \subset \partial K$ such that $\int_E \eta\tilde \varphi\ds = \int_E \tilde \varphi\ds = 0$.

Using similar arguments, we obtain for $\Xi_2$ that
\begin{equation}\label{eq:Xi1}
		\Xi_2  \leq \|\div (\eta \tilde \varphi) \|_{R}\|\xi\|_{R}		
    \leq (1+C_\eta C_\mathrm{P})\|\nabla \tilde \varphi\|_{R}\|\xi\|_{R}.
\end{equation}
To continue the latter estimate, we need a local bound for $\|\xi\|_{K}$. To this end, we test~\cref{eq:LODbasis1} for any element $K$ with $v_\xi \in (H^1_0(K))^n$ chosen such that  $\div v_\xi = \xi$ holds locally in~$K$ and $\|\nabla v_\xi\|_K \leq C_\mathrm{L}'\|\xi\|_K$ is satisfied, cf.~\cref{eq:ladlemma}. This results in
\begin{align}
	\label{eq:bubbleprickforp}
	\|\xi\|_K ^2 = a(\tilde \varphi,v_\xi)\leq C_a C_\mathrm{L}' \|\nabla \tilde \varphi\|_K\|\xi\|_K.
\end{align}
Summing the latter bound over all elements with $K\subset R$, we obtain an estimate for $\|\xi\|_{R}$ which can be inserted in~\cref{eq:Xi1} to conclude
\begin{align*}
	\Xi_2 \leq (1+C_\eta C_\mathrm{P})C_aC_\mathrm{L}'\|\nabla \tilde \varphi\|_{R}^2.
\end{align*}

To estimate $\Xi_3$, we apply the Cauchy--Schwarz inequality for all faces ${E\in\mathcal{F}_H^i (R)}$, where $\mathcal{F}_H^i (R)$ denotes the set of faces inside $R$, which yields that
\begin{equation*}
	 \Xi_3 =
	-\sum_{E \in \mathcal{F}_H^i (R)} \sum_{k = 1}^n | \lambda_{E, k} | \int_E  \eta\tilde{\varphi} \cdot e_k n\sigma  
	\leq \sum_{E \in \mathcal{F}_H^i (R)} \sum_{k = 1}^n | \lambda_{E, k}|  \, |E|_{n-1}^{1/2}  \| \tilde{\varphi} \cdot e_k \|_E,
\end{equation*} 
where $|\cdot|_{n-1}$ denotes the $(n-1)$-dimensional volume. 
We derive a bound for  $| \lambda_{E, k}|$ for all faces $E\in \mathcal{F}_{H}^{i}(R)$ by testing \cref{eq:LODbasis1} with the bubble function $b_{E, k}$, which gives
\begin{equation}
	\label{eq:boundlambda}
	 | \lambda_{E, k} | \leq C_a \| \nabla \tilde{\varphi} \|_{\omega_E} \| \nabla
	b_{E, k} \|_{\omega_E} \leq C_aC_\mathrm{b}H^{- n / 2} \| \nabla \tilde{\varphi} \|_{\omega_E}.
\end{equation}
Here, we used the second estimates of \cref{eq:a} and \cref{eq:bubbleest} and  that $b(b_{E,k},\xi) = 0$ since the divergence of $b_{E,k}$ is piecewise constant and $\xi$ has zero element averages.

Noting that   $|E|_{n-1}  \le H^{n-1}$, the latter estimate can be used to show that
\begin{align*}
	\Xi_3   
	&\leq C_aC_\mathrm{b} H^{-1/2} \sum_{E \in \mathcal{F}_H^i (R)} \sum_{k = 1}^n \| \nabla \tilde{\varphi} \|_{\omega_E}   \| \tilde{\varphi} \cdot e_k \|_E\\
	&	\leq C_a C_\mathrm{b}
	\bigg(n  \sum_{E \in \mathcal{F}_H^i (R)} \| \nabla \tilde{\varphi} \|_{\omega_E}^2 \bigg)^{1/2}  
	\bigg(\frac 1H \sum_{E \in \mathcal{F}_H^i (R)} \| \tilde{\varphi} \|_{E}^2 \bigg)^{1/2},
\end{align*}
where we applied the discrete Cauchy--Schwarz inequality. Using the  trace inequality from \cref{lem:trace}, we obtain that
\begin{equation*}
 \Xi_3  \leq 2C_aC_\mathrm{b} \sqrt{nC_\mathrm{T}(1+C_\mathrm{P})C_\mathrm{P}}\| \nabla \tilde{\varphi} \|_{R}^2.
\end{equation*}

Inserting the above estimates for $\Xi_1$, $\Xi_2$, and $\Xi_3$ into \cref{eq:decaymanipulations} and noting that $R$ can be rewritten as $R = (\Omega\setminus \Nb^{\ell-1}(F)) \setminus (\Omega \setminus \Nb^\ell(F))$, we obtain that
\begin{align*}
	\|\nabla \tilde \varphi\|_{\Omega \setminus \Nb^ \ell(F)}^2\leq C \|\nabla \tilde \varphi\|_{R}^2 =   C \|\nabla \tilde \varphi\|_{\Omega \setminus \Nb^{\ell-1}(F)}^2 - C \|\nabla \tilde \varphi\|_{\Omega\setminus \Nb^ \ell(F)}^2
\end{align*}
with the constant
\begin{equation*}
	C \coloneqq c_a^{-1}\big(C_a'(1 + C_{\eta} C_\mathrm{P}) +(1+C_\eta C_\mathrm{P})C_aC_\mathrm{L}'+2C_aC_\mathrm{b} \sqrt{nC_\mathrm{T}(1+C_\mathrm{P})C_\mathrm{P}}\big)>0.
\end{equation*}
As a direct consequence, we obtain that
\begin{align*}
	\|\nabla \tilde \varphi\|_{\Omega\setminus \Nb^ \ell(F)} \leq \sqrt{\frac{C}{1+C}}\, \|_\nabla \tilde \varphi\|_{\Omega \setminus \Nb^{ \ell-1}(F)},
\end{align*}
which, after iterating, yields that
\begin{equation*}
	\|\nabla \tilde \varphi\|_{\Omega\setminus \Nb^{\ell}(F)} \leq\bigg(\frac{C}{1+C}\bigg)^{\ell/2} \|\nabla \tilde \varphi\|_\Omega= \exp(-c  \ell )\|\nabla \tilde \varphi\|_\Omega
\end{equation*}
for the constant $c \coloneqq \frac{1}{2}\log \tfrac{1+C}{C}>0$. This concludes the proof.
\end{proof}

Motivated by the exponential decay of the prototypical LOD basis functions defined in \cref{eq:LODprotbasis}, we will introduce localizations of them. To formulate the corresponding problems, we introduce localized versions of the spaces $V$ and $Q$ as
\begin{align}
	\label{eq:defVQl}
	\begin{split}
		V_F^\ell &\coloneqq \{v \in V \with \operatorname{supp}(v) \subset \Nb^\ell(F)\},\\
		Q_F^\ell &\coloneqq \{\chi \in Q \with \operatorname{supp}(\chi) \subset \Nb^\ell(F)\}.
	\end{split}
\end{align}
Furthermore, we denote by $R_F^\ell \subset \mathbb R^N$ the subspace consisting of vectors $\mu$ with $\mu_{F,j} = 0$ for all faces not contained in the interior of the patch $\Nb^\ell(F)$, where we recall that $\mu_{F,j}$ is a notation for the entries of the vector $\mu$.

Given an oversampling parameter $\ell$, the localized basis functions, denoted by
\begin{equation*}
	\{\tilde \varphi_{F,j}\with F \in \mathcal F_H^i,\; j = 1,\dots,n\},
\end{equation*}
are determined by the  problems:
Seek $(\tilde \varphi_{F,j}^\ell,\xi_{F,j}^\ell,\lambda) \in V_F^\ell\times Q_F^\ell\times R_F^\ell$  such that
\begin{subequations}
	\label{pbphiEloc} 
	\begin{align}
		&\quad\qquad \qquad a (\tilde \varphi_{F,j}^\ell, v)& +\quad  &b(v,\xi_{F,j}^\ell)          & +\quad  &c(v,\lambda) & =\quad  &0,\qquad \quad\quad &\label{eq:LODlocbasis1}\\
		& \quad\qquad \qquad  b(\tilde \varphi_{F,j}^\ell,\chi)                &   &&   &             & =\quad  &0,\qquad\quad\qquad &\label{eq:LODlocbasis2}\\
		&\quad\qquad\qquad c(\tilde \varphi_{E,j}^\ell,\mu) &   &                  &   &             & =\quad  &\mu_{F, j}\qquad \quad\quad\label{eq:LODlocbasis3}&
	\end{align}
\end{subequations}
holds for all $(v,\chi,\mu) \in V_F^\ell\times Q_F^\ell\times R_F^\ell$. This saddle point problem is well-posed, which can be shown  using similar arguments as in the proof of \cref{le:protbasis}.

In the following, we will frequently use the bound
\begin{equation}
	\label{Boundphil} \| \nabla \tilde{\varphi}_{F, j}^{\ell} \|_{\Omega}
	\leq c_a^{-1}C_a C_\mathrm{b} H^{- n / 2},
\end{equation}
which can be proved by testing equation~\cref{eq:LODlocbasis1} first with $\tilde{\varphi}_{F, j}^{\ell}$ and then with $b_{F, j}$, and combining the resulting estimates.

With the localized basis functions at hand, we are now able to define a localized counterpart of the operator $\mathcal R$ from \cref{eq:defR}. This operator, denoted by $\mathcal R^\ell$, is a projection onto the span of the localized basis functions with the property that it preserves face integrals. More precisely, it is for any $v \in Z$ defined  as 
\begin{equation}
	\label{eq:defRl}
	\mathcal R^\ell v \coloneqq \sum_{F \in \mathcal F_H^i} \sum_{j = 1}^n  
	\int_F v \cdot e_j \ds\, \tilde \varphi_{F,j}^\ell.
\end{equation}
We emphasize that, unlike $\mathcal R$, this projection is not $a$-orthogonal.

The following theorem proves that $\mathcal R^\ell$ approximates $\mathcal R$ exponentially well. Note that using this result for  the bubble function $b_{F,j}$ gives an exponential approximation result for the corresponding localized and prototypical basis functions.

\begin{theorem}[Localization error]\label{thm:approxell}
	There exist a constant $C_\mathrm{R}>0$ independent of~$H$ and $\ell$ such that for all $v \in Z$ and $\ell \in \mathbb N$, it holds that 
	\begin{equation}
		\label{eq:locerr}
		\|\nabla (\mathcal R - \mathcal R^\ell)v\|_\Omega \leq C_\mathrm{R}\ell^{n/2} \exp(-c\ell)\big(\|\nabla v\|_\Omega + H^{-1}\|v\|_\Omega\big).
	\end{equation}
\end{theorem}

\begin{proof}
	We consider an arbitrary but fixed $v \in Z$ and abbreviate $e \coloneqq (\mathcal{R}-\mathcal{R}^{\ell})v$ and  $v_{F, j} : = \int_F v \cdot e_j \ds$, which allows us to write
	\begin{equation}
		\label{eq:decomposition} e = \sum_{F \in \mathcal{F}_H^i} \sum_{j = 1}^n
		v_{F, j} \hspace{0.17em} (\tilde{\varphi}_{F, j} - \tilde{\varphi}_{F,
			j}^{\ell}),
	\end{equation}
	using the definitions of $\mathcal{R}$ and $\mathcal{R}^{\ell}$.
	Applying the first bound from~\cref{eq:a} then gives
	\begin{align}
		c_a  \| \nabla e\|_{\Omega}^2 & \leq \hspace{0.17em} \sum_{F \in
			\mathcal{F}_H^i} \sum_{j = 1}^n v_{F, j} a (\tilde{\varphi}_{F, j} -
		\tilde{\varphi}_{F, j}^{\ell}, e) .  \label{eq:sum}
	\end{align}	
	To estimate the terms on the right-hand side of the latter inequality,
	we fix a face $F \in \mathcal{F}_H^i$ and index $j \in \{1,\dots,n\}$ and recall the definition of the cut-off function in \cref{eq:eta}, which we now denote by  
	$\eta_F$. Note
	that, by the definition of the space~$\tilde{Z}_H$ and since $e \in W$, it
	holds that $a (\tilde{\varphi}_{F, j}, e) = 0$. Using this and observing that $(1 -
	\eta_F) e \in V_F^{\ell}$ is an admissible test function in equation~(\ref{eq:LODlocbasis1}), we obtain that
	\begin{align}\label{eq:sum2}
		a (\tilde{\varphi}_{F, j} - \tilde{\varphi}_{F, j}^{\ell}, e) & = - a
		(\tilde{\varphi}_{F, j}^{\ell}, (1 - \eta_F) e + \eta_F e) 
		\\
		& = b ((1 - \eta_F) e, \xi_{F,j}^\ell) + c ((1 - \eta_F) e, \lambda) - a
		(\tilde{\varphi}_{F, j}^{\ell}, \eta_F e) 
		\notag\\
		& \eqqcolon \Xi_1 + \Xi_2 + \Xi_3 . 
		\notag
	\end{align}
	
In the following, we estimate the terms $\Xi_1$, $\Xi_2$, and $\Xi_3$ separately. To estimate~$\Xi_1$, we note that the divergence of $(1 - \eta_F)e$ is piecewise constant on $\Omega \setminus R_F$ and that~$\xi_{F,j}^\ell$ has zero element averages. This together with~\cref{eq:boundeta,lem:Poincare}~gives
	\begin{align}
		\Xi_1 \leq \| \div (1 - \eta_F) e\|_{R_F} \| \xi_{F,j}^\ell \|_{R_F} \leq (1 +
		C_{\eta_F} C_{\mathrm{P}})  \| \nabla e\|_{R_F} \| \xi_{F,j}^\ell \|_{R_F} . 
		\label{eq:Xi11}
	\end{align}
	To derive a $L^2$-bound for $\xi_{F,j}^\ell$, we proceed similarly as in~\cref{eq:bubbleprickforp}, which yields that
	\begin{align*}
		\| \xi_{F,j}^\ell \|_{R_F} \leq C_a C_{\mathrm{L}}'  \| \nabla \tilde{\varphi}_{F,
			j}^{\ell} \|_{R_F},
	\end{align*}
	and inserting this into~\cref{eq:Xi11} gives
	\begin{align*}
		\Xi_1 \leq C_a C_{\mathrm{L}}'  (1 + C_{\eta_F} C_{\mathrm{P}})  \| \nabla
		e\|_{R_F}  \| \nabla \tilde{\varphi}_{F, j}^{\ell} \|_{R_F} .
	\end{align*}
		
	To estimate $\Xi_2$, we apply the discrete Chauchy--Schwarz
	inequality, the trace inequality from \cref{lem:trace}, the bound $|E|_{n - 1} \leq H^{n - 1}$, and the local Poincaré-type inequality from 
	\cref{lem:Poincare}, which yields that
	\begin{align*}
		\Xi_2 & \leq \bigg( \sum_{E \subset R_F} \sum_{k =
			1}^n \lambda_{E, k}^2 \bigg)^{1 / 2}  \bigg(  \sum_{E \subset R_F} |E|_{n
			- 1} \|e\|_E^2 \bigg)^{1 / 2}\\
		& \leq \bigg( \sum_{E \subset R_F} \sum_{k = 1}^n \lambda_{E, k}^2
		\bigg)^{1 / 2}  \bigg( C_{\mathrm{T}}  \sum_{E \subset R_F} |E|_{n - 1}
		\|e\|_K (\| \nabla e\|_K + H^{- 1} \|e\|_K) \bigg)^{1 / 2}\\
		& \leq \bigg( \sum_{E \subset R_F} \sum_{k = 1}^n \lambda_{E, k}^2
		\bigg)^{1 / 2} \sqrt{nC_{\mathrm{T}} C_{\mathrm{P}}  (1 +
			C_{\mathrm{P}})} H^{n / 2}  \| \nabla e\|_{R_F} .
	\end{align*}
	To derive a bound for $| \lambda_{E, k} |$, we proceed as in \cref{eq:boundlambda}, but locally, which results in
	\begin{align*}
		| \lambda_{E, k} | \leq |a (\tilde{\varphi}_{F, j}^\ell, b_{E, k}) | \leq C_a
		C_{\mathrm{b}} H^{- n / 2}  \| \nabla \tilde{\varphi}_{F, j}^\ell \|_{\omega_E}
		.
	\end{align*}
	Using this estimate, the above estimate for  $\Xi_2$ can be continued as 
	\begin{align*}
		\Xi_2 \leq C_a C_{\mathrm{b}} n \sqrt{C_{\mathrm{T}} C_{\mathrm{P}}  (1 +
			C_{\mathrm{P}})}  \| \nabla e\|_{R_F}  \| \nabla
		\tilde{\varphi}_{F, j}^{\ell} \|_{R_F} .
	\end{align*}
	
	To estimate $\Xi_3$, we again use similar arguments as in \cref{eq:corrected}, which yields that	\begin{align*}
		\Xi_3 \leq C_a'  \| \nabla (\eta_F e)\|_{R_F}  \| \nabla
		\tilde{\varphi}_{F, j}^{\ell} \|_{R_F} \leq C_a'  (1 + C_{\eta_F}
		C_{\mathrm{P}})  \| \nabla e\|_{R_F}  \| \nabla \tilde{\varphi}_{F,
			j}^{\ell} \|_{R_F}
	\end{align*} 
		with $C_a' \coloneqq (\nu_\mathrm{max} + C_\mathrm{P}^2\sigma_\mathrm{max})$. 
	
	Inserting the above estimates for $\Xi_1$, $\Xi_2$, and $\Xi_3$ into \cref{eq:sum2} gives the bound
	\begin{equation*}
		 a (\tilde{\varphi}_{F, j} - \tilde{\varphi}_{F, j}^{\ell}, e) \leq C
		\| \nabla e\|_{R_F}  \| \nabla \tilde{\varphi}_{F, j}^{\ell} \|_{R_F}
	\end{equation*}
	with
	\begin{equation*}
		C \coloneqq \big( (C_a C_{\mathrm{L}}' + C_a')  (1 + C_{\eta_F} C_{\mathrm{P}}) +C_a C_{\mathrm{b}} n \sqrt{C_{\mathrm{T}} C_{\mathrm{P}}  (1 +
			C_{\mathrm{P}})} \big).
	\end{equation*}
	We can now apply the exponential decay result of \cref{thm:dec} to $\tilde{\varphi}_{F, j}^{\ell}$ instead of~$\tilde \varphi_{F,j}$, where we replace~$\Omega$ in the statement of the theorem by $\Nb^{\ell} (F)$. Using this and bound \cref{Boundphil} for the localized basis function then gives
	\begin{align*}
		a (\tilde{\varphi}_{F, j} - \tilde{\varphi}_{F,
			j}^{\ell}, e) &\leq C \exp (- c \ell)  \| \nabla
		e\|_{R_F}  \| \nabla \tilde{\varphi}_{F, j}^{\ell} \|_{\Nb^{\ell} (F)}
		\leq C' H^{- n / 2} \exp (- c \ell)  \| \nabla e\|_{R_F},
	\end{align*}
	where we abbreviated  $C' \coloneqq c_a^{-1}C_aC_\mathrm{b}C$. 
	
	It remains to sum the above estimate over all faces $F \in \mathcal F_H^i$ and indices $j \in \{1,\dots,n\}$ as in \cref{eq:sum}. Using the discrete Cauchy--Schwarz inequality, the trace inequality from \cref{lem:trace}, and Young's inequality, we get that
	\begin{align*}
		c_a  \| \nabla e\|_{\Omega}^2 & \leq \sqrt{n} C' \exp (- c
		\ell) H^{- n / 2} \bigg( \sum_{F \in \mathcal{F}_H^i} | F |_{n - 1}
		\|v\|_F^2 \bigg)^{1 / 2} \bigg( \sum_{F \in \mathcal{F}_H^i} \| \nabla
		e\|_{R_F}^2 \bigg)^{1 / 2} \\
		& \leq \sqrt{nC_{\mathrm{T}}} C' C_{\mathrm{ol}} \ell^{n/ 2}
		 \exp (- c \ell)  \| \nabla e\|_{\Omega} (\| \nabla
		v\|_{\Omega} + H^{- 1} \|v\|_{\Omega}),
	\end{align*}
	where $C_{\mathrm{ol}}>0$ is a constant only depending on the
	regularity of the mesh $\mathcal{T}_H$. This proves the desired estimate with the constant $C_\mathrm{R} \coloneqq a_a^{-1}\sqrt{nC_{\mathrm{T}}} C' C_{\mathrm{ol}}$.
\end{proof}

\section{Localized multiscale method}
\label{sec:pracmethod}

In this section, we introduce the proposed multiscale method for heterogeneous Stokes problems. Its approximation space, denoted by $\tilde Z_H^\ell$, is defined as the span of the localized basis functions~\cref{pbphiEloc}, i.e.,
\begin{equation}
	\tilde Z_H^\ell \coloneqq \operatorname{span}\{\tilde \varphi_{F,j}^\ell\with F \in \mathcal F_H^i,\; j = 1,\dots,n\}.
\end{equation}
The proposed multiscale method then seeks $(\tilde u_H^\ell,\tilde p_H^\ell) \in \tilde Z_H^\ell\times M_H$ such that
	\begin{subequations}
	\label{LODid2loc} 
	\begin{align}
		\qquad \qquad \qquad &a(\tilde u_H^\ell,\tilde v_H^\ell)& +\quad &b(\tilde v_H^\ell,\tilde p_H^\ell) &=\quad &(f,\tilde v_H^\ell)_\Omega,&\qquad \qquad \qquad\quad\label{eq:LODdivfreeloc1}\\
		\qquad \qquad \qquad&b(\tilde u_H^\ell,\tilde q_H^\ell)&   &       &=\quad &0&\qquad \qquad \qquad\label{eq:LODdivfreeloc2}
	\end{align}
\end{subequations}
 holds for all $(\tilde v_H^\ell,\tilde q_H^\ell) \in (\tilde Z_H,M_H)$.

The remainder of this section is devoted to the error analysis of this method. We emphasize that the pressure approximation $\tilde p_H^\ell$ is piecewise constant, and therefore, e.g., first-order convergence can only be expected if $p\in H^1(\Omega)$. However, such regularity requirements are generally not satisfied for heterogeneous Stokes problems  (or, if they are, the $H^1$-norm of $p$ can be very large). As a remedy, we introduce a post-processing step that uses the local pressure contributions~$\xi_{F,j}^{\ell}$ computed together with the LOD basis functions in~\cref{pbphiEloc}, as
\begin{equation}
	\label{plpp} \tilde{p}_H^{\ell, \text{pp}} \coloneqq \tilde{p}_H^{\ell} + \sum_{F
		\in \mathcal{F}_H^i} \sum_{j = 1}^n c_{F,j}  \hspace{0.17em} \xi_{F,
		j}^{\ell},
\end{equation}
where $c_{F,j}$ is the coefficient of $\tilde \varphi_{F,j}^\ell$ in the representation of $\tilde{u}_H^{\ell}$.
The following theorem proves the well-posedness of the method~\cref{LODid2loc} and its uniform convergence properties for the velocity and post-processed pressure approximations under minimal regularity assumptions, provided that the $\ell$ is chosen sufficiently large. In addition, it is proved that the piecewise constant (not post-processed) pressure approximation $\tilde p_H^\ell$ converges exponentially 
to $\Pi_H p$ as $\ell$ is increased.

\begin{theorem}	[Localized method]\label{thm:errestloc}
	The localized multiscale method \cref{LODid2loc} is well-posed. Furthermore, there exist constants $C_u, C_u', C_p, C_p'> 0$ independent of $H$ and~$\ell$ such that for any right-hand side $f \in H^s(\Omega)$ with $s \in \{0,1\}$, it holds that
	\begin{align}
		\| \nabla (u - \tilde{u}_H^{\ell})\|_{\Omega} & \leq C_u  \big( H^{1+s}|f|_{s,\Omega} + H^{- 1} \ell^{n/ 2} \exp (- c \ell) \|f\|_{\Omega}\big),  \label{eq:errestpracH1}\\
		\|u - \tilde{u}_H^{\ell} \|_{\Omega} & \leq C_u'  \big(H  + H^{- 1} \ell^{n/ 2} \exp (- c \ell)\big)\| \nabla (u - \tilde{u}_H^{\ell})\|_{\Omega}  \label{eq:errestpracL2}
		\end{align}
		for the  velocity approximation,
		where we recall that $|\cdot|_{s,\Omega}$ denotes the $H^s(\Omega)$-seminorm. For the pressure approximation, we have the error estimates
				\begin{align}
			\| \Pi_H p - \tilde{p}_H^{\ell} \|_{\Omega} &\leq C  H^{- 1} \ell^{n/ 2} \exp (- c \ell) (1 + H^{- 1} \ell^{n/ 2} \exp (- c \ell)) \|f\|_{\Omega},
		\label{eq:errestpracp}\\
		\|p - \tilde{p}_H^{\ell, \mathrm{pp}} \|_{\Omega} &\leq C_p  \big(H + H^{- 1}
		\ell^{n/ 2} \exp (- c \ell)\big) \big( 1  + H^{- 1}
		\ell^{n / 2} \exp (- c \ell) \big) \|f\|_{\Omega}.
		\label{eq:errestpracpp}
	\end{align}

\end{theorem}

\begin{proof}
	We begin this proof by showing the inf--sup condition
	\begin{equation}
		\label{infsuplH} \adjustlimits \inf_{\tilde{q}_H^{\ell} \in M_H}
		\sup_{\tilde{v}_H^{\ell} \in \tilde{Z}_H^{\ell}}  \frac{|b
			(\tilde{v}_H^{\ell}, \tilde{q}_H^{\ell}) |}{\| \nabla \tilde{v}_H^{\ell}
			\|_{\Omega} \| \tilde{q}_H^{\ell} \|_{\Omega}} \geq \tilde{c}_b^{\ell}
	\end{equation}
	for some constant $\tilde{c}_b^{\ell} > 0$ to be specified later, which
	implies the well-posedness of problem~\cref{LODid2loc}. To this end, we first show the continuity of the operator~$\mathcal{R}^{\ell}$. Using the continuity of~$\mathcal{R}$, cf.~(\ref{eq:Rcont}), the approximation result from \cref{thm:approxell}, and the Poincaré--Friedrichs inequality on $\Omega$, it follows that
	\begin{align*}
		\| \nabla \mathcal{R}^{\ell} v\|_{\Omega} & \leq \big( \sqrt{C_a / c_a} + C_\mathrm{R} \ell^{n/ 2} \exp
		(- c \ell) (1 + C_{\mathrm{PF}} H^{- 1}) \big)  \| \nabla v\|_{\Omega}.
	\end{align*}	
	Given any $\tilde{q}_H^{\ell} \in M_H$, we denote by $v \in V$ the
	function satisfying $\div v = \tilde{q}_H^{\ell}$ and $\| \nabla
	v \|_{\Omega} \leq C_\mathrm{L} \| \tilde{q}_H^{\ell} \|_{\Omega}$,
	cf.~(\ref{eq:ladlemmaOmega}). Choosing $\tilde{v}_H^{\ell} \coloneqq
	\mathcal{R}^{\ell}  v$, we obtain noting that \cref{eq:identityIH} also holds for the operator  $\mathcal R^\ell$ that $b( \tilde{v}_H^{\ell},\tilde q_H^\ell) =
b(v,\tilde q_H^\ell) = \|\tilde q_H^\ell\|_\Omega^2$. The desired inf--sup condition~\cref{infsuplH} then follows with the constant 
	\begin{equation}
		\label{ctillb} \tilde{c}_b^{\ell} = C_\mathrm{L}^{- 1} \big( \sqrt{C_a
			/ c_a} + C_\mathrm{R} \ell^{n/ 2} \exp (- c \ell) (1 +
		C_{\mathrm{PF}} H^{- 1}) \big)^{- 1}.
	\end{equation}

	To prove the $H^1$-error estimate for the velocity approximation, we denote by $\tilde{Z}_H^{\ell, 0} \coloneqq \{ v \in \tilde Z_H^\ell \with \div v = 0\}$ the  subspace of divergence-free functions of the localized approximation space. Problem~\cref{LODid2loc} can then be equivalently reformulated as the unique solution to: seek $\tilde{u}_H^{\ell}	\in \tilde{Z}_H^{\ell,0}$ such that
	\begin{equation*}
		 a (\tilde{u}_H^{\ell}, \tilde{v}_H^{\ell}) = (f,
		\tilde{v}_H^{\ell})_{\Omega}
	\end{equation*}
	holds for $ \tilde{v}_H^{\ell} \in
	\tilde{Z}_H^{\ell, 0}$.
	
	 Similarly, also the prototypical multiscale method~\cref{LODid2} can be equivalently reformulated in the corresponding subspace of divergence-free functions defined as $\tilde{Z}_H^0 \coloneqq \{v \in \tilde{Z}_H \with \div v = 0\}$, i.e., we seek $\tilde{u}_H \in \tilde Z_H^0$ such that 
	 \begin{equation*}
	 	 a (\tilde{u}_H, \tilde{v}_H) = (f, \tilde{v}_H)_{\Omega}  
	 \end{equation*}holds for all $
	 \tilde{v}_H \in \tilde{Z}_H^{0}$.
	   
	   Interpreting $\tilde{u}_H^{\ell} \in \tilde{Z}_H^{\ell, 0}$ as a non-conforming, non-consistent approximation of $\tilde{u}_H \in \tilde{Z}_H^{ 0}$, we can apply Strang's second lemma, cf.~\cite[Lem.~2.25]{ErG04}, which gives
\begin{align}\label{eq:strang2}
	\begin{split}
  \| \nabla (\tilde{u}_H - \tilde{u}_H^{\ell})\|_{\Omega} \leq \big( 1 +
  c_a^{-1}C_a \big) &\inf_{\tilde{v}_H^{\ell} \in \tilde{Z}_H^{\ell,
   0}} \| \nabla (\tilde{u}_H - \tilde{v}_H^{\ell})\|_{\Omega} 
  \\
  + c_a^{- 1}&  \sup_{\tilde{v}_H^{\ell} \in \tilde{Z}_H^{\ell,
   0}}  \frac{| (f, \tilde{v}_H^{\ell})_{\Omega} - a (\tilde{u}_H,
  \tilde{v}_H^{\ell}) |}{\| \nabla \tilde{v}_H^{\ell} \|_{\Omega}}.
\end{split}
\end{align}
To estimate the infimum on the right-hand side of \cref{eq:strang2}, we choose $\tilde{v}_H^{\ell} =\mathcal{R}^{\ell} u$, and note that $\mathcal{R}^{\ell} u \in \tilde{Z}_H^{\ell, 0}$. The latter holds since property \cref{eq:identityIH} also holds for operator~$\mathcal R^\ell$. Furthermore, since we can identify $\tilde{u}_H = \mathcal R u$, cf.~\cref{thm:convergenceprot}, we have that $\tilde{u}_H - \tilde{v}_H^{\ell} = (\mathcal{R}-\mathcal{R}^{\ell}) u$, which allows to apply \cref{thm:approxell}. This gives
\begin{equation*}
 \inf_{\tilde{v}_H^{\ell} \in \tilde{Z}_H^{\ell,  0}} \| \nabla
(\tilde{u}_H - \tilde{v}_H^{\ell})\|_{\Omega} \leq C_\mathrm{R}
\ell^{n/ 2} \exp (- c \ell) (\| \nabla u\|_{\Omega} + H^{- 1}
\|u\|_{\Omega}).
\end{equation*}

The supremum on the right-hand hand side of \cref{eq:strang2} can be estimated noting that for any $\tilde{v}_H^{\ell} \in \tilde{Z}_H^{\ell, 0}$ and $\tilde{v}_H \in \tilde{Z}_H^{ 0}$ it holds that
\begin{equation*}
	 | (f, \tilde{v}_H^{\ell})_{\Omega} - a (\tilde{u}_H, \tilde{v}_H^{\ell}) |
	= | (f, \tilde{v}_H^{\ell} - \tilde{v}_H)_{\Omega} - a (\tilde{u}_H,
	\tilde{v}_H^{\ell} - \tilde{v}_H) |.
\end{equation*}
Given $\tilde{v}_H^{\ell} \in \tilde{Z}_H^{\ell, 0}$ we choose $\tilde{v}_H =\mathcal{R} \tilde v_H^\ell$ and use that $\tilde{v}_H - \tilde{v}_H^{\ell} = (\mathcal{R}-\mathcal{R}^{\ell}) \tilde v_H^\ell \in W$, which implies that $a (\tilde{u}_H, \tilde{v}_H^{\ell} - \tilde{v}_H)=0$. Applying \cref{thm:approxell} and the Poincaré--Friedrichs inequality then yields  the estimate
\begin{multline*}
   \sup_{\tilde{v}_H^{\ell} \in \tilde{Z}_H^{\ell,  0}} 
    \frac{| (f,
   \tilde{v}_H^{\ell})_{\Omega} - a (\tilde{u}_H, \tilde{v}_H^{\ell})|}
   {\| \nabla \tilde{v}_H^\ell \|_{\Omega}}
   \leq C_\mathrm{PF}   C_\mathrm{R}
   \ell^{n/ 2} \exp (- c \ell) (1 + H^{- 1} C_{\text{PF}})\| f  \|_{\Omega}.
\end{multline*}

We are now ready to continue~\cref{eq:strang2}. Applying the Poincaré--Friedrichs inequality again and using the stability estimate $\| \nabla u\|_\Omega \leq c_a^{- 1} C_{\mathrm{PF}} \|f\|_{\Omega}$, we obtain that
\begin{equation}
  \label{Cea} \| \nabla (\tilde{u}_H - \tilde{u}_H^{\ell})\|_{\Omega}
  \leq a_a^{-1}\big( 2 + c_a^{-1}{C_a} \big)  C_\mathrm{R} (C_\mathrm{PF} + H^{- 1} C_{\mathrm{PF}}^2) \ell^{n
  / 2} \exp (- c \ell)  \| f
  \|_{\Omega}.
\end{equation}
Estimate~\cref{eq:errestpracH1} can now be concluded with the constant $C_u \coloneqq \max\{C_\mathrm{e}',c_a^{-1}(2+c_a^{-1}C_a)C_\mathrm{R}C_\mathrm{PF}(1+C_\mathrm{PF})\}>0$ using the triangle inequality and the $H^1$-convergence result for the prototypical method from \cref{thm:convergenceprot}. Note that in order to simplify the constant, we have used that $\Omega$ is of unit size, which means that $H \leq 1$.

Using similar arguments as above, one can also show that the piecewise constant pressure approximation $\tilde p_H^\ell$ converges exponentially to $\Pi_H p$ in the $L^2$-norm as $\ell$ is increased. More specifically, considering the combined problem obtained by adding up equations \cref{eq:LODdivfreeloc1,eq:LODdivfreeloc2}, one can apply Strang's second lemma in a similar way as above and receive, among other things, estimate \cref{eq:errestpracp}.
The proof of this estimate will not be given in detail here, because it is similar to the proof above.

To prove the $L^2$-error estimate for the velocity approximation, we use an Aubin--Nietsche-type duality argument. Denoting by $(\bar{u}, \bar{p} ) \in V \times M $ the solution to~\cref{eq:weakstokes} for the right-hand side $g \coloneqq u - \tilde{u}_H^{\ell}$, we obtain for any $v \in Z_H^{\ell, 0}$ that
	\begin{align*}
		\|u - \tilde{u}_H^{\ell} \|_{\Omega}^2 & = (g, u -
		\tilde{u}_H^{\ell})_{\Omega}= a (u - \tilde{u}_H^{\ell}, \bar{u} - v) \leq C_a
		\| \nabla (u - \tilde{u}_H^{\ell}) \|_{\Omega} \| \nabla (\bar{u} -
		v) \|_\Omega.
	\end{align*} 
	Choosing $v$ as the approximation of the proposed  method to 	$\bar{u}$, and applying the already established $H^1$-error estimate 	\cref{eq:errestpracH1} for $s = 0$ to $\bar{u}$, we get that
	\begin{equation*}
		 \|u - \tilde{u}_H^{\ell} \|_{\Omega}^2 \leq C_a C_u  \big( H + H^{- 1} \ell^{n/ 2} \exp (- c \ell) \big)\| \nabla (u - \tilde{u}_H^{\ell}) \|_{\Omega}\|g\|_{\Omega}.
	\end{equation*}
	Estimate~\cref{eq:errestpracL2} can then be concluded with the constant $C_u' \coloneqq C_a C_u$.
	
	To prove the $L^2$-error estimate for the post-processed pressure approximation, we consider an arbitrary but fixed element $K \in \TH$ and a function $v \in (H^1_0 (K))^n$ to be specified later. For any $F \in \mathcal{F}_H^i$ and $j \in \{1, \dots, n\}$, we then test equation~\cref{eq:LODlocbasis1} determining the respective localized basis functions with $v$. Multiplying the resulting equation by $c_{F,j}$, the coefficient of $\tilde \varphi_{F,j}^\ell$ in the basis representation of~$\tilde u_H^\ell$, cf.~\cref{plpp}, and summing up yields that
	\begin{equation*}
		 a_K (\tilde{u}_H^{\ell}, v) + b_K  (v, \tilde{p}_H^{\ell, \text{pp}} -
		\tilde{p}_H^{\ell}) = 0,
	\end{equation*}
	where the subscript  denotes the restriction of bilinear forms $a$ and $b$ to~$K$. 

Moreover, testing \cref{eq:weakstokes} with the same $v$ and using that by the divergence theorem there holds  $b_K(v,\Pi_H p) = (\Pi_Hp)|_K \int_K \div v\dx = 0$, we obtain that
\begin{equation*}
	 a_K (u, v) + b_K  (v, p - \Pi_H p) = (f, v)_K.
\end{equation*}
Subtracting the latter two equations results in
\begin{equation}
	\label{eq:test}
	 a_K (u - \tilde{u}_H^{\ell}, v) + b_K  (v, p - \Pi_H p -
	(\tilde{p}_H^{\ell, \text{pp}} - \tilde{p}_H^{\ell})) = (f, v)_K.
\end{equation}
Abbreviating $q \coloneqq p -\Pi_H p - (\tilde{p}_H^{\ell, \text{pp}} - \tilde{p}_H^{\ell})$, the Ladyzhenskaya lemma, cf.~\cref{eq:ladlemma}, asserts the existence of a function $v \in (H^1_0 (K))^n$ such that $\div v = q$ holds locally in~$K$ and which satisfies $\|\nabla v \|_K \leq C_L^\prime \|q\|_K$. Choosing this~$v$ in equation~\cref{eq:test} and using the uniform coefficient bounds \cref{eq:unifboundcoeff} as well as the local Poincaré--Friedrichs inequality $\| v \|_K \leq H \| \nabla v \|_K$ for all $v \in (H^1_0(\Omega))^n$, we obtain  that
	\begin{align*} 
		\|q
	\|_K^2 &= (f, v)_K - a_K (u - \tilde{u}_H^{\ell}, v)\\
	&\leq C_\mathrm{L}' \big(H \| f
	\|_K + \sigma_\mathrm{max}H \|u - \tilde{u}_H^{\ell}\|_K + \nu_\mathrm{max}\| \nabla (u - \tilde{u}_H^{\ell}) \|_K \big)\| q\|_K.
	\end{align*}
	Summing this inequality over all $K \in \TH$, we arrive at
	\begin{equation}
		\label{eq:q}
	 \| q
	\|_{\Omega} \leq \sqrt{3}C_\mathrm{L}' \big( H \| f \|_{\Omega} + \sigma_\mathrm{max}H\|u-\tilde u_H^\ell\|_\Omega  +\nu_\mathrm{max}\|
	\nabla (u - \tilde{u}_H^{\ell}) \|_{\Omega}\big).
	\end{equation}
	The triangle  inequality inequality finally yields that
	\begin{equation*}
		 \| p - \tilde{p}_H^{\ell, \text{pp}} \|_{\Omega} \leq \| \Pi_H p -
		\tilde{p}_H^{\ell} \|_{\Omega} + \|q\|_{\Omega},
	\end{equation*}
	where the first term can be bounded using \cref{eq:errestpracp}
	and the second term using \cref{eq:q} and estimates \cref{eq:errestpracH1,eq:errestpracL2} for $s = 0$. This proves estimate \cref{eq:errestpracpp} with the constant $C_p \coloneqq C+ \sqrt{3}C_\mathrm{L}'(1+\sigma_\mathrm{max}C_u' + \nu_\mathrm{max}C_u)$, which concludes the proof.
\end{proof}

If the oversampling parameter is suitably coupled to the coarse mesh size, one obtains uniform spatial convergence, as can be seen in the following corollary. 

\begin{corollary}[Uniform convergence]
	Let $f \in H^s(\Omega)$. Then, if the oversampling parameter $\ell$ is increased logarithmically with the coarse mesh size $H$, i.e., $\ell \in \mathcal O(\log 1/H)$, we have for a constant $C>0$ independent of $H$ and $\ell$ that
	\begin{align*}
		\| \nabla (u - \tilde{u}_H^{\ell})\|_{\Omega} \leq CH^{1+s},\qquad \|u - \tilde{u}_H^{\ell} \|_{\Omega} \leq C H^{2+s},\qquad \|p - \tilde{p}_H^{\ell, \mathrm{pp}} \|_{\Omega} \leq C H.
	\end{align*}
\end{corollary}

\section{Implementation and numerical experiments}
\label{sec:numexp}

In this section, we discuss the implementation of the proposed multiscale method and present numerical experiments that support the theoretical results of this paper.

\subsection*{Implementation}
For a practical implementation of the method, the local but still infinite-dimensional patch problems \cref{pbphiEloc} have to be discretized. For this purpose, we use a fine-scale discretization based on the Crouzeix--Raviart FEM (CR-FEM), cf.~\cite{Crouzeix1973}. The CR-FEM is particularly well suited as a fine-scale discretization since its piecewise constant pressure approximation space allows  definition~\cref{eq:defZ}, the resulting reformulated Stokes problem~\cref{eq:reformulation}, and the problems defining the localized basis functions~\cref{pbphiEloc} to be easily adapted to the fully discrete setting. The fully discrete velocity approximation obtained by the~\cref{LODid2loc}, is divergence-free in the weak discrete sense, i.e., its divergence vanishes on all elements of the fine mesh. This does not mean, however, that it is divergence-free on $\Omega$, since the CR space is non-conforming.  For simplicity, we assume that the global fine mesh, denoted by $\mathcal T_h$, is obtained by (multiple) uniform red refinement of the coarse mesh~$\TH$. Note that the fine mesh needs to resolve all microscopic features of the coefficients to obtain a reliable approximation. The patch problems \cref{pbphiEloc} are then discretized using discrete versions of the infinite-dimensional spaces $V_F^\ell$ and $Q_F^\ell$, cf.~\cref{eq:defVQl}, defined on the local meshes obtained by restricting $\mathcal T_h$ to the respective patches (the global mesh $\mathcal{T}_h$ is never used directly). To practically enforce that the Lagrange multipliers~$\xi_{F,j}^\ell$ in \cref{pbphiEloc} have zero element averages, another Lagrange multiplier, which is piecewise constant with respect to $\TH$, must be added. 

 The approximation results for the prototypical multiscale method from  \cref{thm:convergenceprot} can be  transferred to the case of a CR fine-scale discretization with some minor modifications. For example, one needs to generalize the Poincaré-like inequality from \cref{lem:Poincare} to the case of fine-scale CR functions which are not $H^1_0(\Omega)$-conforming. This can be easily done using the concept of conforming companions to CR-functions, cf.~\cite[Chap.~5.2]{Gallistl2014Adaptive}. With this tool at hand, one can redo the proof of \cref{thm:convergenceprot}, now estimating the error between the fully discrete prototypical LOD solution and the fine-scale CR-FEM solution. The only part of the proof that deserves special attention is the integration by parts in~\cref{trickyCR}. Here, due to the non-conformity, we need to apply integration by parts on each element of the fine mesh $\mathcal{T}_h$, resulting in jump terms at the boundaries of the elements. These jump terms can be estimated by $h\|f\|_K\|e\|_K$, which needs to be added to the right-hand side of~\cref{eq:est2}. The exponential decay and approximation result from \cref{thm:dec,thm:approxell} can be easily adapted to the fully discrete setting by inserting an appropriate interpolation to the CR space on $\mathcal{T}_h$ wherever necessary, cf.~\cite[Chap.~4.4]{MalP20}. Adapting the proof of \cref{thm:errestloc} using the aforementioned results in the fully discrete setting, then gives, for example, the $H^1$-error estimate:
 \begin{equation}
	\label{eq:fsdisc}
	\|\nabla (u_h-\tilde u_{H,h}^\ell)\|_\Omega \leq  C\big(H^2 \|\nabla f\|_\Omega  + (H^{-1}\ell^{n/2}\exp(-c\ell) + h)\|f\|_\Omega\big),
\end{equation}
where $\tilde u_{H,h}^\ell$ and $u_h$ denote the fully discrete LOD solution and fine-scale CR-FEM solution, respectively.  and  $C>0$ is a constant independent of $H$, $\ell$, and $h$. An error estimate against the continuous solution can be inferred using the triangle inequality, estimate \cref{eq:fsdisc}, and classical a priori convergence results for the CR-FEM. Note that in error estimates against the continuous solution, the term $C h\|f\|_\Omega$ is dominated by the error of the fine-scale CR-FEM. Nevertheless, some of our numerical experiments where we compute the error against the fine-scale CR-FEM solution (not detailed here), confirm the presence of this term in the error estimates. 

\subsection*{Numerical experiments}
\begin{figure}
		\centering
		\begin{subfigure}[t]{0.5\textwidth}
			\centering
			\includegraphics[height=.89\linewidth]{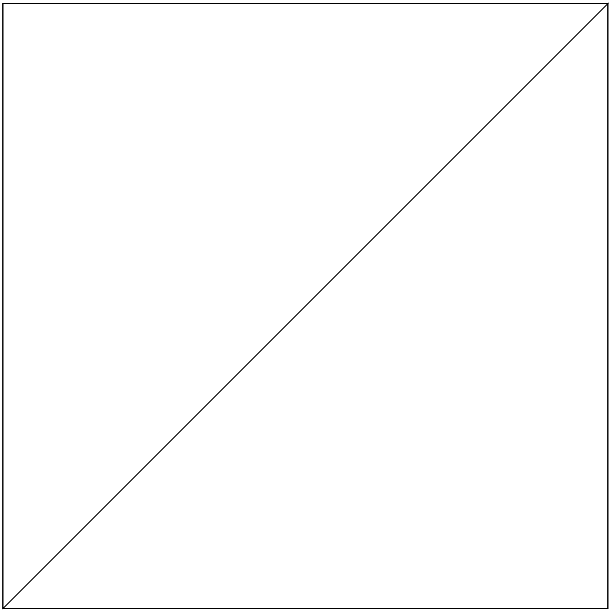}
		\end{subfigure}%
		~ 
		\begin{subfigure}[t]{0.5\textwidth}
			\centering
			\includegraphics[height=.9\linewidth]{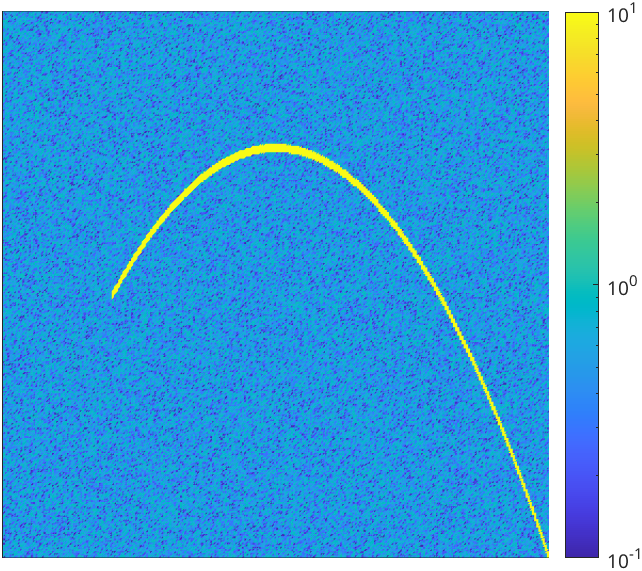}
		\end{subfigure}
	\caption{Initial mesh for the mesh generation (left). Multiscale coefficient used in all numerical experiments (right).}
	\label{fig:coeff}
\end{figure}

We consider the domain $\Omega = (0,1)^2$ and introduce a hierarchy of meshes generated by uniform red refinement of the initial mesh shown in \cref{fig:coeff} (left). For simplicity, we denote the meshes in the hierarchy by $\mathcal T_{2^0}, \mathcal T_{2^{-1}}, \dots$, where the subscript refers to the side length of the squares formed by joining opposing triangles. The coefficient $\nu$ is chosen to be piecewise constant with respect to the mesh~$\mathcal T_\epsilon$  with element values obtained as realizations of independent random variables uniformly distributed in the interval $[0.1, 1]$. Note that~$\epsilon$ is  assumed to be a negative power of two. For elements whose midpoints have a distance less than $4\epsilon$ from a parabola, the corresponding element values are set to~10; see \cref{fig:coeff} (right). The coefficient $\sigma$ is set to zero for simplicity.

Note that all numerical experiments presented below can be reproduced using the code available at  \url{https://github.com/moimmahauck/Stokes_LOD_CR}. 

\subsubsection*{Exponential decay of basis functions}
For the first numerical experiment we use the coefficient $\nu$ as described above for the value $\epsilon = 2^{-6}$. The mesh for the fine-scale discretization is chosen to be $\mathcal T_{2^{-8}}$, which sufficiently resolves the coefficient. This relatively large fine mesh size is necessary to compute the prototypical LOD basis functions needed to evaluate the localization errors.

\begin{figure}
	\centering
	\begin{subfigure}[t]{0.5\textwidth}
		\centering
		\includegraphics[height=.9\linewidth]{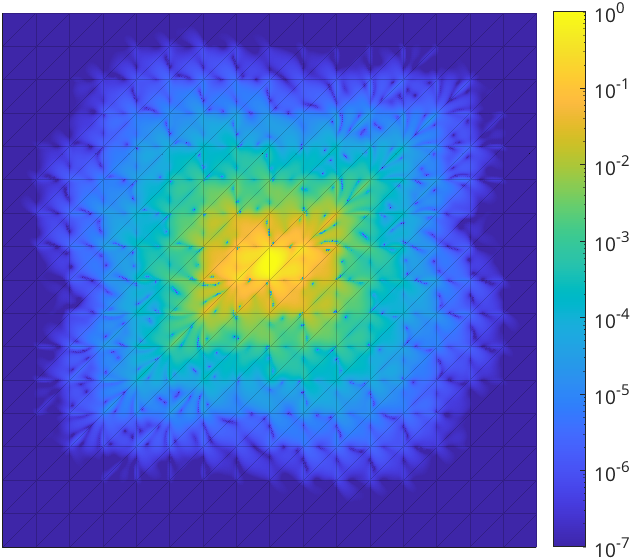}
	\end{subfigure}%
	~ 
	\begin{subfigure}[t]{0.5\textwidth}
		\centering
		\includegraphics[height=.9\linewidth]{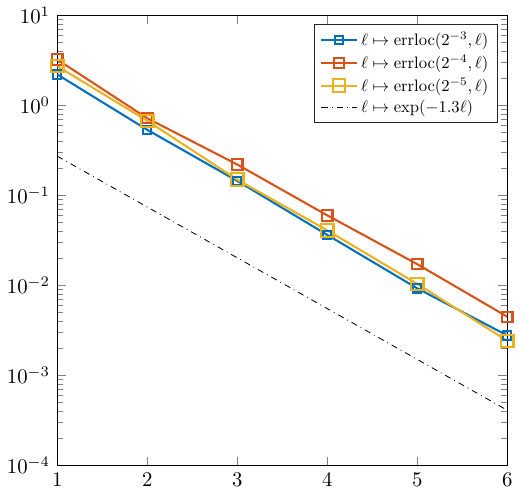}
	\end{subfigure}
	\caption{Decay of the modulus of a prototypical LOD basis function, plotted using a logarithmic color scale (left). $H^1$-errors of the localized approximation of the prototypical LOD basis functions for several coarse mesh sizes $H$, plotted as a function of the oversampling parameter $\ell$ (right).}
\label{fig:localization}
\end{figure}

In \cref{fig:localization} (left), we illustrate the modulus of an exemplary basis function using a logarithmic color scale. A trained eye observes an exponential decay of the modulus with respect to the underlying coarse mesh, which we indicated in light gray. This supports the exponential decay result of \cref{thm:dec}. Next, we numerically investigate the localization error when replacing a prototypical basis function by its localized counterpart. For a given coarse mesh size $H$ and localization parameter $\ell$, we define the $H^1$-localization error as
\begin{equation*}
	\mathrm{errloc}(H,\ell) \coloneqq \max_{F \in \mathcal F_H^i}\max_{j = 1,\dots,n} \|\nabla (\tilde \varphi_{F,j,h}-\tilde \varphi_{F,j,h}^\ell)\|_\Omega.
\end{equation*}
In \cref{fig:localization} (right), one clearly observes an exponential decay of the $H^1$-norm localization error as the localization parameter $\ell$ is increased. This supports the exponential approximation result from~\cref{thm:approxell}. 

\subsubsection*{Optimal order convergence}

For the second numerical experiment, we use the coefficient $\nu$ from above for the value ${\epsilon = 2^{-8}}$ and choose the mesh $\mathcal T_{2^{-10}}$ for the fine-scale discretization. Furthermore, as right-hand side we use the function
\begin{equation*}
	f(x,y) \coloneqq  (-y,x)^T.
\end{equation*}

\begin{figure}
		\centering
	\begin{subfigure}[t]{0.5\textwidth}
		\centering
		\includegraphics[height=.9\linewidth]{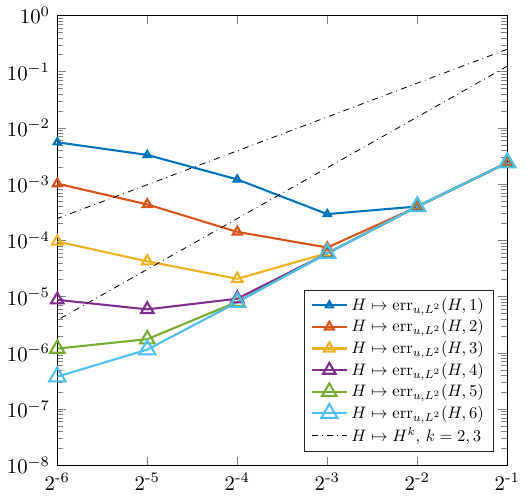}
	\end{subfigure}%
	~ 
	\begin{subfigure}[t]{0.5\textwidth}
		\centering
		\includegraphics[height=.9\linewidth]{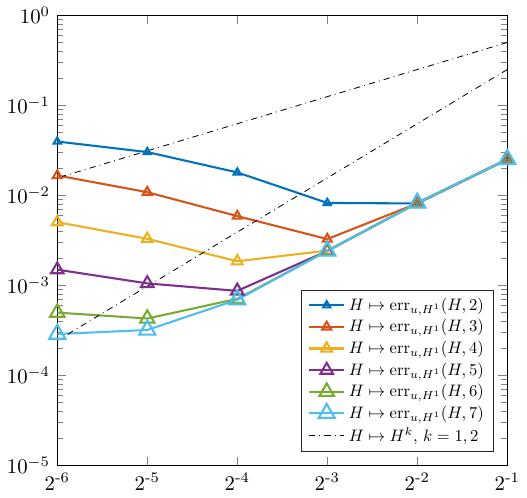}
	\end{subfigure}
	\caption{$L^2$-errors (left) and $H^1$-errors (right) for the velocity approximation for several oversampling parameters $\ell$, plotted as functions of the coarse mesh size $H$.}
	\label{fig:convu}
\end{figure}

In the following, we investigate the  errors
\begin{align*}
	\mathrm{err}_{u,H^1}(H,\ell) &\coloneqq \|\nabla(u_h - \tilde u_{H,h}^\ell)\|_\Omega, & \quad 
	\mathrm{err}_{u,L^2}(H,\ell)&\coloneqq \|u_h - \tilde u_{H,h}^\ell\|_\Omega,\\
		\mathrm{err}_{p,L^2}(H,\ell)&\coloneqq \|p_h - \tilde p_{H,h}^{\ell,\mathrm{pp}}\|_\Omega, & \quad
			\mathrm{err}_{\Pi_H p,L^2}(H,\ell)&\coloneqq \|\Pi_H p_h - \tilde p_{H,h}^\ell\|_\Omega,
\end{align*}
where we recall that above $(u_h,p_h)$ denotes the reference solution computed on the fine mesh. 
For the $L^2$- and $H^1$-errors of the velocity approximation, we observe in \cref{fig:convu} the (almost) third and second order convergence, respectively, provided that the localization parameter is chosen sufficiently large. Recalling that ${f \in H^1(\Omega), }$, this is in line with the prediction from \cref{thm:errestloc}.  For a fixed oversampling parameter, we observe that after a certain error level is reached, the error increases again as the mesh size is decreased. This is a well-known effect that occurs for some LOD methods such as \cite{MalP14,Maier2021}. It can be overcome with a more sophisticated localization strategy; see, e.g., \cite{HenP13,Hauck2022,Dong2023}. Such an improved localization strategy will be investigated in future~work. 

\begin{figure}
			\centering
	\begin{subfigure}[t]{0.5\textwidth}
		\centering
		\includegraphics[height=.9\linewidth]{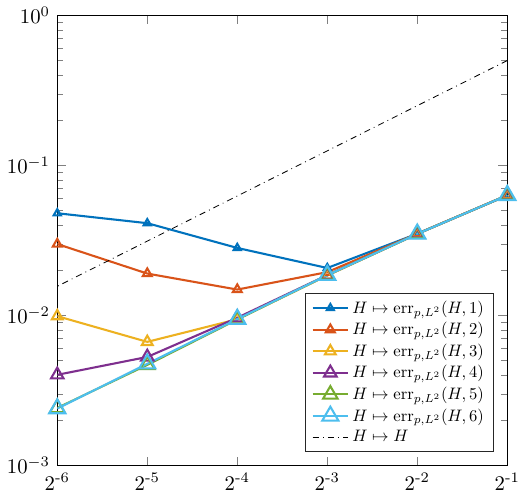}
	\end{subfigure}%
	~ 
	\begin{subfigure}[t]{0.5\textwidth}
		\centering
		\includegraphics[height=.9\linewidth]{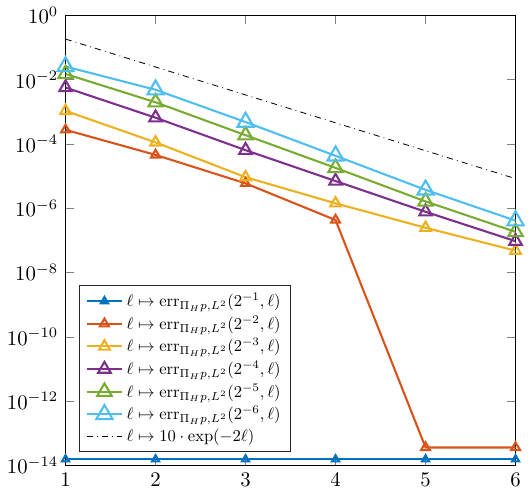}
	\end{subfigure}
	\caption{$L^2$-errors of the post-processed pressure approximation for several oversampling parameters $\ell$, plotted as functions of the coarse mesh size $H$ (left). $L^2$-errors of the pressure approximation computed with respect to $\Pi_H p$ for several coarse mesh sizes~$H$, plotted as a function of the oversampling parameter $\ell$ (right).}
	\label{fig:convp}
\end{figure}

Now we turn to the pressure approximation. For the $L^2$-error of the post-processed pressure approximation we observe in \cref{fig:convp} (left) the first-order convergence, again provided that the localization parameter is chosen sufficiently large. This observation is in line with \cref{thm:errestloc}. For the piecewise constant pressure approximation $\tilde p_{H,h}^\ell$, we observe the exponential $L^2$-convergence towards~$\Pi_H p_h$ in \cref{fig:convp} (right), which is also consistent with \cref{thm:errestloc}. Note that the outliers are due to combinations of $H$ and $\ell$ where all patches are global, which implies that the pressure approximation coincides with $\Pi_H p_h$, cf. \cref{thm:convergenceprot}.

\section*{Acknowledgment}
The authors would like to thank Axel M{\aa}lqvist for helpful discussions on the construction and analysis of the proposed method.

\appendix
\section{Collection of frequently used bounds}
\begin{lemma}[Local Poincaré-type inequality]\label{lem:Poincare}
	There exists $C_\mathrm{P}>0$ independent of $H$ such that, for all $K \in \TH$ and all $v \in H^1(K)$ satisfying $\tint_F v \ds = 0$ for at least one face with $F \subset \partial K$, it holds that
	\begin{equation}
		\label{eq:Poincarélocal}
		\|v\|_K \leq C_\mathrm{P} H \|\nabla v\|_K.
	\end{equation}
\end{lemma}
\begin{proof}
	This result can be derived from \cite[Lem.~B.66]{ErG04} using a transformation to the reference element and the corresponding estimates in~\cite[Lem.~1.101]{ErG04}.
\end{proof}

\begin{lemma}[Trace inequality]\label{lem:trace}
	There exists a constant $C_\mathrm{T}>0$ independent of~$H$ such that, for all $K \in \TH$ and $v \in H^1(K)$, it holds for any face with $F \subset \partial K$ that
	\begin{equation*}
		\|v\|_F^2 \leq C_\mathrm{T}\big(\|\nabla v\|_K + H^{-1}\|v\|_K\big)\|v\|_K.
	\end{equation*} 
\end{lemma}
\begin{proof}
	The proof of this result can be done as in \cite[Lem.~1.49]{PiE12}, invoking the quasi-uniformity we assumed for the sequence of meshes considered.
\end{proof}

\bibliographystyle{alpha}
\bibliography{bib}
\end{document}